\documentclass[letterpaper, 12pt]{amsart}

\usepackage[T1]{fontenc}

\usepackage[margin=1in]{geometry}

\usepackage{amscd}

\usepackage{mathrsfs}

\usepackage{amsmath}

\usepackage{amsfonts}

\usepackage{amssymb}

\usepackage{amsthm}

\usepackage{arydshln}

\usepackage{fancyhdr}

\usepackage{verbatim}

\usepackage{tikz-cd}

\usepackage[all]{xy}

\usepackage{makecell}

\usepackage{ulem}
\usepackage{cancel}

\usepackage{color}

\usepackage[colorlinks=true, linkcolor=blue, citecolor=blue,
pagebackref=true]{hyperref}

\usepackage{mathtools} \usepackage{etoolbox}
\patchcmd{\section}{\scshape}{\bfseries}{}{} \makeatletter
\renewcommand{\@secnumfont}{\bfseries} \makeatother

\theoremstyle{definition}

\theoremstyle{plain} \newtheorem{theorem}{Theorem}[section]

\theoremstyle{plain} \newtheorem{lemma}[theorem]{Lemma}

\theoremstyle{plain} \newtheorem{proposition}[theorem]{Proposition}

\theoremstyle{plain} 

\theoremstyle{plain} \newtheorem{corollary}[theorem]{Corollary}

\theoremstyle{remark} \newtheorem*{remark}{Remark}

\theoremstyle{definition} \newtheorem{definition}[theorem]{Definition}

\theoremstyle{definition} \newtheorem*{definition*}{Definition}

\theoremstyle{definition} 

\theoremstyle{remark} 

\makeatletter \renewenvironment{proof}[1][\proofname]
{\par\pushQED{\qed}\normalfont\topsep6\p@\@plus6\p@\relax\trivlist\item[\hskip\labelsep\bfseries#1\@addpunct{.}]\ignorespaces}{\popQED\endtrivlist\@endpefalse}
\makeatother

\newcommand{\con}{\textrm{con}}

\newcommand{\RR}{\mathbb{R}}

\newcommand{\Pfull}{P_{\textrm{full}}}
\newcommand{\Pmod}{P_{\textrm{mod}}}
\newcommand{\Sfull}{S_{\textrm{full}}}

\providecommand{\Pfulln}[1]{S_{\textrm{full}}(\Sigma^R, \tilde\psi, #1)}
\providecommand{\Pmodn}[1]{S_{\textrm{mod}}(\Sigma, \psi, #1)}
\providecommand{\Pordn}[1]{S_{\textrm{full}}(\Sigma, \psi, #1)}

\renewcommand{\leq}{\leqslant} \renewcommand{\geq}{\geqslant}

\providecommand{\ceil}[1]{\lceil#1\rceil}

\DeclarePairedDelimiter{\floor}{\lfloor}{\rfloor}

\newcommand{\bi}{{\bf i}}
\newcommand{\ri}{{\rm i}}

\newcommand{\rk}{{\rm k}}
\newcommand{\rr}{{\rm r}}

\newcommand{\brr}{{\bf r}}

\newcommand{\rmq}{{\rm q}}
\newcommand{\ba}{{\bf a}}
\newcommand{\rb}{{\rm b}}
\newcommand{\ra}{{\rm a}}
\newcommand{\rc}{{\rm c}}

\newcommand{\rw}{{\rm w}}
\newcommand{\bj}{{\bf j}}
\newcommand{\rj}{{\rm j}}

\newcommand{\ul}{\underline{\ell}}
\renewcommand{\l}{{\partial}}
\newcommand{\Lmin}{L_{\min}}

\newcommand{\R}{\mathbb R}
\newcommand{\N}{\mathbb N}

\def\1int{{[0,1]}}

\title[Shrinking targets in affine iterated function systems]{Shrinking targets in affine iterated function systems }

\author[Koivusalo, Liao, Morris, Rams]{Henna Koivusalo, Lingmin Liao, Ian D. Morris and Micha{\l} Rams}


\address{Fry Building, Woodland Rd, Bristol BS8 1UG, UK}
\email{henna.koivusalo@bristol.ac.uk}

\address{School of Mathematics and Statistics, Wuhan University, Wuhan, Hubei, 430072, China}
\email{lmliao@whu.edu.cn}

\address{School of Mathematical Sciences, Queen Mary University of London, Mile End Road, London E1 4NS, UK}
\email{i.morris@qmul.ac.uk}

\address{Institute of Mathematics, Polish Academy of Sciences, ul. \'Sniadeckich 8, 00-656 Warszawa,  Poland} 
\email{rams@impan.pl}

\date{\today}

\dedicatory{}

\begin{document}


\begin{abstract}
We calculate the Hausdorff dimension of shrinking target sets of symbolic balls in generic affine iterated function systems with matrix norms bounded by $\tfrac 12$. Relative to earlier works of e.g. Barany and Troscheit, and Koivusalo and Ramirez, we impose only very mild condition, known as complete reducibility, on the iterated function system. We also generalise these earlier works in that the size of the target in the current work is allowed to depend on the trajectory hitting it. The formula for Hausdorff dimension we obtain is presented as a zero point of certain pressure function. 
\end{abstract}


\maketitle
\thispagestyle{empty}
\section{Background}\label{Sec:backg}

Questions about infinite recurrence and recurrence rates are classical in the study of dynamical systems. The shrinking target sets consist of those points which have infinitely many returns to a fixed target sequence. More precisely, given a dynamical system $(X,T)$, and a sequence of targets $(U_n)$, the {\it shrinking target set } is 
\[
\mathcal S(T, (U_n))=\{x\in X\mid T^n(x)\in U_n \textrm{ for infinitely many }n\}. 
\]
The term `shrinking targets' was coined by Hill and Velani in the 90's \cite{HV95, HV97}, although the concept arises naturally e.g. in the context of dynamics and Diophantine approximation. In \cite{HV95} Hill and Velani studied, in particular, expanding maps $T$ on the Riemann sphere, with $X$ the corresponding Julia set, and the targets $U_n=B(z_0, r_n)$ chosen to be balls shrinking at some rate around a fixed centre point. In this set-up they proved that the Hausdorff dimension of $\mathcal S(T, (U_n))$ is given by a zero of a suitable pressure function. 

Just two years later, in \cite{HV97} Hill and Velani introduced a variant of this set-up, whereby the target sequence $U_n$ depends on the orbit of $x$. More precisely, in order for $x\in \mathcal S(T, (U_n))$, they required infinitely many $n$, $y\in T^{-n}(z_0)$, for which
\[
T^n(x)\in B(y, r_n(x)). 
\]
Here the sequence $(r_n(x))$ is chosen according to a H\"older continuous potential on $X$. To make a distinction to the first case, where the target sequence was fixed for all $x$, we call this second variant, with $r_n(x)$ dependent on the orbit point, a {\it path-dependent shrinking target set}. Also in the path-dependent set-up, the Hausdorff dimension of $\mathcal S(T, (U_n))$ is given by a pressure function. These works have since then been followed by a rich literature of shrinking target results, for example for $\beta$-transformations\cite{BW14}, for the Gauss map \cite{LWWX14}, and for countable Markov maps \cite{R11}. 

In this article we focus on shrinking target sets in the context of {\it iterated function systems} (IFSs). An iterated function system is a finite collection of contractive maps $\{f_1, \dots, f_N\}$ on $\R^d$, which by a classical result of Hutchinson \cite{H81} always yields a unique, non-empty, compact invariant set $\Lambda=\cup_{i=1}^N f_i(\Lambda)$. The geometric properties of $\Lambda$ depend delicately on the regularity and scaling properties of $f_i$. The particular class of IFSs that will be of interest to us in this work, are affine IFSs, of the form $x\mapsto f_i(x)=T_i(x)+v_i$, with $T_i$ linear and $v_i\in \R^d$ a translation. There are two dynamical systems that arise naturally in this context: first of all, for many self-affine sets there is an expanding map $E$ on $\Lambda$, with $f_i$ as its local inverses. The second dynamical system is the left shift map on the coding space $\Sigma=\{1, \dots, N\}^{\N}$ which can be used to code the points of $\Lambda$. Much of the  interesting theory for the shrinking target problem in the self-affine situation arises from the interplay between the geometry of the set $\Lambda$ (or, indeed, the maps of the generating iterated function system) and the dynamics on the coding space. For more detail, see Section \ref{sec:preliminaries}. Versions of the shrinking target problem for the non-path-dependent case for iterated function systems have been covered in \cite{KR18, BT, JK, BR, U02}. 

In this paper, we study the path-dependent shrinking target sets in affine iterated function systems, proving as Theorem \ref{main} that their Hausdorff dimension is the zero of a pressure function. The applicability of our work to the path-dependent set-up is one of the major novelties of our work. However, there is another one, as our proof significantly expands the class of affine iterated function systems, for which the Hausdorff dimension of the shrinking target set can be computed. 

Let us give some context, to understand the fundamental difference between our work and those in the past literature. A major dynamical difficulty in the dimension theory of non-conformal iterated function systems in general, and the study of the shrinking target sets in particular, is that the operation of combining the contractions in affine iterated function systems is usually only sub-multiplicative and not multiplicative. Natural covers of the attractor using cylinders in the shift space may therefore change in an abrupt and somewhat irregular manner between different scales. In the literature there are several workarounds, and in particular a lot of modern theory (e.g. \cite{BT, KR18}) relies on various weak quasi-multiplicativity conditions, which can be shown to be satisfied under certain irreducibility assumptions. These conditions hold generically for choices of linear maps $T_1, \dots, T_N$, and so for large classes of affine IFSs. 

However, we are able to substantially overcome this limitation, working in the context of invertible affine iterated function systems whose linear parts only need to satisfy an even weaker {\it complete reducibility} condition. Completely reducible systems are precisely those systems which are simultaneously block diagonalisable with all diagonal blocks being irreducible. In particular (perhaps counter-intuitively given the terminology), all irreducible systems are completely reducible, since these systems correspond to the case of a single irreducible diagonal block. At the other extreme, all simultaneously diagonalisable systems are also completely reducible. In particular, our work is a significant improvement to the current literature. For definitions, see Section \ref{sec:preliminaries}. 

The novel technique we use should be applicable more broadly in affine thermodynamic formalism, not just for the shrinking problem, which is why we have separated the key ideas into the relatively self-contained Section \ref{sec:modular}. Namely, we write the shift space modularly in multiplicative blocks, and demonstrate that the  impact of this difference in expression to e.g. pressure is negligible. This idea was inspired by the proofs in \cite{BT}. In Section \ref{sec:targetpressure} this modular space idea is combined with the recently discovered fact (see \cite{MorKae}) that via a change of basis, every affine IFS with complete reducibility can be approximately decomposed into a collection of quasi-multiplicative parts. In particular, in this case we find a measure on the shift space which is well-distributed. The only other assumption we have on the affine iterated function system, that the norms of the linear parts are bounded from above by $\tfrac 12$, is inherited from a transversality techniques of \cite{F88}, which comes into play when the symbolic problem is interpreted geometrically on $\Lambda$. This is the reason for the result being formulated for generic translations of the affine maps.

\section{Preliminaries and the statement of results}\label{sec:preliminaries}
\subsection{Symbolic space}
Denote $A=\{1, \dots, N\}$, $\Sigma:= A^\mathbb{N}$ and let $\sigma$ be the left shift operator on $\Sigma$. Then the pair $(\Sigma, \sigma)$ is a dynamical system called the full shift dynamics on the alphabet $A$ of $N$ symbols. Let $\Sigma_n:= A^n$ be the set of words of length $n$. We denote infinite words by bold letters $\bi, \, \bj,  \, \ba$ and so on, and finite words by $\ri, \, \rj, \, \ra$ and so on. The set $\Sigma_*:= \cup_{n=1}^\infty \Sigma_n$ is the collection of all finite words. 
For any $\ri=i_1\cdots i_n\in \Sigma_n$, denote by $[\ri]$ the cylinder corresponding to $\ri$,  i.e.,
\[
[\ri]:=\big\{\bj=(j_1, j_2, \dots)\in \Sigma: \ j_1=i_1, \dots, j_n=i_n\big\}.
\]
The length $n$ of $\ri \in \Sigma_n$ is denoted by $|\ri|$. For $\bi\in \Sigma$, denote by $\bi|_n:=i_1i_2\cdots i_n$ the finite word composed of the first $n$ symbols of $\bi$, and by $\bi|^m_n:= i_{m+1}i_{m+2}\cdots i_{n}$ the finite word composed of the symbols between the positions $m$ and $n$. Such a finite word $\bi|^m_n$ is called a subword of $\bi$. For convenience, for a positive real number $\ell$, we write $\bi|_\ell$ for $\bi|_{\lfloor \ell\rfloor}$, where $\lfloor \cdot \rfloor$ denotes the integer part. One can easily interpret the similar symbols $\ri|_n$, $\ri|^m_n$ and $\ri|_{\lfloor \ell\rfloor}$ for a finite word $\ri\in \Sigma_*$.

\subsection{Symbolic shrinking targets}
We will investigate a variation of the shrinking target problem for $(\Sigma, \sigma)$. A {\it length sequence function} is a function $\ell: \Sigma \rightarrow (\R^+)^\N$ defined as $\bi \mapsto \ell(\bi)=(\ell_n(\bi))_{n\geq 1}$. The value $\ell(\bi)$, which depends on an infinite word $\bi$, will give the lengths (sizes) of the shrinking targets. The center of the shrinking targets will be an infinite word $\bj\in \Sigma$. Then for any $\bj, \bi\in \Sigma$, we define a sequence of finite words, i.e., a family of targets: $$\bj(\bi, \ell):=\big(\bj|_{\ell_n(\bi)}\big)_{n\geq 1}.$$   
For a length sequence function $\ell$ and $\bj\in \Sigma$, define the following {\it symbolic path-dependent shrinking target set}
\begin{align}\label{def-STS}
R(\bj, \ell)=\Big\{\bi\in \Sigma\mid \sigma^n(\bi)\in [\bj|_{\ell_n(\bi)}]\textrm{ for infinitely many }n\Big\}. 
\end{align}

\subsection{Iterated function systems and shrinking targets}\label{Subsec:2.3}

Let $\{f_1, \dots, f_N\}$ be a collection of affine contractions, that is, let $T_1, \dots, T_N$ be linear contractions and $v_1, \dots, v_N\in \R^d$, and let $f_i(x)=T_i(x)+v_i$ for $i=1, \dots, N$. This is called an {\it affine iterated function system}. For $\ri=i_1i_2\cdots i_n\in \Sigma_*$, we denote
\[
T_\ri:=T_{i_1}\circ T_{i_2} \circ \cdots \circ T_{i_n}, 
\]
and similarly for $f_\ri$. Recall that by a classical theorem of Hutchinson \cite{H81}, an iterated function system defines a unique, non-empty, compact, invariant self-affine set $\Lambda$ such that 
\[
\Lambda=\bigcup_{i=1}^Nf_i(\Lambda). 
\]
Define the coding map 
\[
\pi: \Sigma\to \Lambda, \ \pi(\bi)=\lim_n f_{\bi|_n}(0), 
\]
and then set the geometric shrinking target set on $\Lambda$ to be
\begin{equation}\label{def-STS-real}
R^*(\bj, \ell):=\pi(R(\bj, \ell)). 
\end{equation}

We say that the affine iterated function system satisfies the strong separation condition if the images $f_i(\Lambda)$ are disjoint. In particular, then the mapping $\pi$ is a bijection, and each $x\in \Lambda$ corresponds to exactly one infinite sequence $\bi=\pi^{-1}(x)=x_1x_2\dots$. The symbols $x_k$ are called the digits of $x$. Further, in this case there is an everywhere well-defined expanding map $E$ on $\Lambda$ with $f_i$ being its local inverses, given by $E: \Lambda \to \Lambda$, $E(x)=f_{x_1}^{-1}(x)$, where $x_1$ is the first digit of $\pi^{-1}(x)$. The map $E$ is conjugate to the associated symbolic shift dynamics. That is we have the following commutative diagram:
 \begin{center}
 \setlength{\unitlength}{1mm}
 \begin{picture}(40,35)
 \put(4,4){$\Lambda$}
 \put(10,5){$\vector(1,0){15}$}
 \put(15,0.7){$E$}
 \put(28,4){$\Lambda$}
 \put(4,26){$\Sigma$}
 \put(6,23){$\vector(0,-1){15}$}
 \put(2.4,15){$\pi$}
 \put(10,27){$\vector(1,0){15}$}
 \put(15,28){$\sigma$}
 \put(28,26){$\Sigma$}
 \put(30,23){$\vector(0,-1){15}$}
  \put(31,15){$\pi$}
 \end{picture}
 \end{center}
 
\vskip 5pt

Hence, the symbolic shrinking target sets $R(\bj, \ell)$ defined above have geometric interpretation as subsets of $\Lambda$, and
\begin{align}
R^*(\bj, \ell)=\pi(R(\bj, \ell))=\Big\{x\in \Lambda\mid E^n(x)\in \pi\big[\bj|_{\ell_n(\pi^{-1}(x))}\big]\textrm{ for infinitely many }n\Big\}. 
\end{align}
If the strong separation condition does not hold, such geometric interpretations cannot be made and, since we aim to prove dimension results for generic translations, imposing a strong separation condition would be unnatural. The target sets $\pi[\bj|_{\ell_n(\pi^{-1}(x)}]$ also make the best sense symbolically. For these reasons, we focus on the set $\pi(R(\bj, \ell))$.

\subsection{Additive and sub-additive potentials}

A {\it potential} is a function $\phi:\Sigma\to \RR$. Together with a potential one considers its Birkhoff sums, for $\bi \in \Sigma$,
\[
S_n\phi(\bi) = \sum_{i=0}^{n-1} \phi(\sigma^i \bi).
\]

A special class of potentials are {\it piecewise constant } potentials whose values depend only on the first symbol, $\phi(\bi)=\phi(i_1)$. For piecewise constant potentials, their Birkhoff sums also depend only on the first finitely many symbols: for $\bi\in \Sigma$,
\[
S_n \phi(\bi) = S_n \phi(i_1,\ldots,i_n).
\]
A {\it sub-additive potential} is a family of functions $\phi_n:\Sigma\to\RR; \ n=1,2,\ldots$ satisfying for $\bi\in \Sigma$
\begin{equation}\label{eqn:subpot}
\phi_{m+n}(\bi) \leq \phi_m(\bi) + \phi_n(\sigma^m\bi).
\end{equation}
Clearly, for any piecewise constant potential $\phi$, setting $\phi_n = S_n\phi$, the inequality \eqref{eqn:subpot} is automatically satisfied and is an equality, hence the notion of sub-additive potential is indeed a generalization of the notion of a potential or, more precisely, a generalization of its family of Birkhoff sums.

The piecewise constant potentials have their analogue among sub-additive potentials: potentials such that for every $n$, $\phi_n$ depends only on the first $n$ symbols. In this case, \eqref{eqn:subpot} takes the following 'concatenating' form: for an element $\bi\in \Sigma$, $\bi=(i_1, i_2, i_3,\dots)$,
\[
\phi_{m+n}(i_1,\ldots,i_{m+n}) \leq \phi_m(i_1,\ldots,i_m) + \phi_n(i_{m+1},\ldots,i_{m+n}).
\]
Note that any function $\psi : \Sigma_* \rightarrow \mathbb{R}$ induces a family of functions $\phi_n$ defined on $\Sigma$, and depending only on the first $n$ symbols in the following natural way:
\[
\phi_n(\bi)=\psi(\ri), \quad \forall \ri \in \Sigma_*, \bi\in [\ri].
\]
Thus, we usually study functions defined on $\Sigma_*$, and call a potential $\psi : \Sigma_* \rightarrow \mathbb{R}$ sub-additive if the induced family $\phi_n$ is sub-additive.

We will also need to discuss {\it sub-multiplicative potentials}. We call a potential $\psi$ on $\Sigma^*$ sub-multiplicative  if it takes positive real values and satisfies 
\[
\psi(\ri\rj)\le \psi(\ri)\psi(\rj). 
\]
Using sub-additive or sub-multiplicative potentials is a matter of convention, as any sub-multiplicative potential can be turned into a sub-additive one by taking a logarithm.

\subsection{Weakly quasi-additive singular value potentials}
To a linear map $T: \R^d\to \R^d$ and a parameter $s>0$ we can associate a singular value function 
\[
\phi^s(T)=\alpha_1\alpha_2\cdots\alpha_{\lfloor s \rfloor}\alpha_{\lfloor s \rfloor+1}^{s-\lfloor s \rfloor},
\]
where $\alpha_i$ are the singular values of $T$  in descending order. The singular value function was first introduced in affine dimension theory by Falconer \cite{F88}. For a fixed collection $T_1, \dots, T_N$  of linear maps on $\R^d$, we write $\phi^s(\ri) $ instead of $\phi^s(T_\ri)$. 

The potential $\phi^s$ on $\Sigma_*$ is sub-multiplicative, and hence $\log \phi^s$ is sub-additive. When investigating affine iterated function systems, one needs to work with the thermodynamic formalism for $\log \phi^s$. Notice that, this potential is defined on some matrix cocycle, generated by contracting maps, which means that $\log \phi^s$ is strictly negative for $s>0$.

Furthermore, the singular value potential quite often has better properties than mere sub-additivity. In particular,  we will be concerned with whether this potential $\psi=\log \phi^s$ is actually {\it weakly quasi-additive}: 
\begin{definition}[weak quasi-additivity]
A potential $\psi$ on $\Sigma_*$ is {\it weakly quasi-additive} when there exist constants $Q,K$ such that for any two finite words $\ri,\rj\in\Sigma_*$ we have some $\rk\in \Sigma_*, |\rk|\leq K$ such that
\begin{equation}\label{eq:weakquasi}
\psi(\ri)+ \psi(\rj)\ge \psi(\ri \rk \rj) \geq Q + \psi(\ri) + \psi(\rj).
\end{equation}
We denote the connecting word corresponding to $\ri, \ra\in \Sigma^*$ by $\con(\ri,\ra)$. If there are more than one such word, we use the one which is smallest in the lexicographical order. We call the potential {\it quasi-additive} if $\con(\ri,\ra)$ can always be taken to be the empty word.
\end{definition}
This property and its stronger variants have been a crucial assumption in much of the past work on the shrinking target problem on self-affine sets \cite{KR18, BT}. 

\begin{remark}
We emphasise at this point that the property of weak quasi-additivity is in fact not a property of $\psi$ alone, but of $\psi$ and the matrix cocycle generated by $T_1, \dots, T_N$ together. However,  since in the context of this article an underlying predefined collection of linear maps is considered fixed, we use slightly imprecise language and call the potential $\psi$ weakly quasi-additive. For example, the singular value potential $\phi^s$ is weakly quasi-additive for matrix cocycles all of whose exterior powers are strongly irreducible, but can fail to be weakly quasi-additive when a single exterior power becomes reducible (see \cite{KM18} and \cite{MS19} respectively).
\end{remark}

The weak quasi-additivity condition is noticeably weaker than quasi-additivity. In \cite{BT} B\'ar\'any and Troscheit proposed a very interesting approach to handling thermodynamic formalism for weakly quasi-additive potentials, which in some sense allows us to reduce the weakly quasi-additive situation to quasi-additive one. In the arguments below we are building on their proof idea. 

\subsection{Statement of the main theorem}

We are now ready to start formulating our main theorem, which concerns the Hausdorff dimension of a {\it path-dependent shrinking target set}. The motivation for setting up the problem in this way is from Hill and Velani \cite{HV97}, Bugeaud and Wang \cite{BW14}, Li, Wang, Wu, and Xu \cite{LWWX14}, Reeve \cite{R11}, et al,  who have treated in the conformal setting the problem of path-dependent shrinking targets, where the target balls are given by a H\"older continuous potential.

We need some assumptions on the length sequence function $\ell$.
\begin{definition}[Assumptions on $\ell$]\label{def:path}
Let $\ell:\Sigma\to (\R^+)^\N$ be a function defined as $$\ell(\bi)=(\ell_n(\bi))_{n\geq 1}, \quad \forall \bi\in \Sigma.$$ 
Assume that $\ell_n(\bi)$ depends only on the first $n$ symbols of $\bi$. Then we can define an associated function $\underline{\ell}:\Sigma_*\to (\R^+)$ by for each $\ri\in \Sigma_*$ setting $\underline{\ell}(\ri) := \ell_{|\ri|}(\ba)$ for any $\ba\in [\ri]$. 
We assume that $\underline{\ell}$ is {\it approximately additive} on finite words: there exists a constant $\kappa>0$ such that for any $\ri, \ra\in \Sigma_*$, 
\begin{align}\label{estimate-ell}
|\underline{\ell}(\ri\ra)-\underline{\ell}(\ri)-\underline{\ell}(\ra)|\leq \kappa. 
\end{align}
We also assume that 
\begin{align}\label{assmp:ln-infinity}
\ell_n(\bi)\to\infty \quad \text{for every}  \ \bi\in\Sigma.
\end{align}
\end{definition}

\begin{remark}
 We can easily find a function $\ell$ satisfying the assumptions in Definition \ref{def:path}. Let $\psi$ be a potential on $\Sigma$ which depends only on the first symbol, and let 
\[
\ell_n(\bi)=\sum_{i=0}^{n-1}\psi\circ \sigma^i(\bi). 
\]
Then  $\ell_n$ only depends on the first $n$ symbols and $\underline{\ell}$ is approximately additive on finite words, with $\kappa=0$. 
\end{remark}

Over the course of this work, we will define many different notions of pressure, with the following notion of shrinking target pressure being one of the most crucial ones. For a glossary of different pressures, see Section \ref{sec:glossary}. 

\begin{definition}\label{def:pressure}
For $\bj\in \Sigma$, and a length sequence function $\ell$, define the following limsup pressure function
\begin{equation}\label{eq:pressure}
P^*(s, \bj)=\limsup_n\tfrac 1n\log \sum_{|\ri|=n}\phi^s(\ri\bj|_{\underline{\ell}(\ri)}). 
\end{equation}
When the limit exists, we denote the pressure by $P(s, \bj)$. 
\end{definition}
Recall here that when $\ul(\ri)$ is not an integer, we interpret $\bj|_{\ul(\ri)}=\bj|_{\lfloor \ul(\ri)\rfloor}$. Also, note that the notation $\ri\bj|_{\underline{\ell}}$ means the word $\ri$ followed by the first $\ell$ symbols of the word $\bj$, it is not the same as $(\ri\bj)|_{\underline{\ell}}$.

We require one further notion in order to state our main theorem. We call a tuple of linear maps $T_1,\ldots,T_N \colon \R^d \to \R^d$ {\it irreducible} if the only linear subspaces $V$ of $\R^d$ which satisfy $T_iV \subseteq V$ for all $i=1,\ldots,N$ are $\{0\}$ and $\R^d$. Following a terminological convention from representation theory, we call $T_1,\ldots,T_N$ {\it completely reducible} if for some invertible linear map $X \colon \R^d \to \R^d$ we may write
\[X T_iX^{-1} = \begin{pmatrix}
T_i^1&0&\cdots&0\\
0&T_i^2&\cdots &0 \\ 
\vdots & \vdots& \ddots&\vdots\\
0&0&\cdots &T^k_i
\end{pmatrix}
\]
for all $i=1,\ldots,N$, where for each $j=1,\ldots,k$ the tuple $T_1^j,\ldots,T_N^j$ is an irreducible tuple of square matrices of some dimension $d_j$, say. We observe that a tuple of diagonal matrices is completely reducible since one may take $k=d$; at the other extreme,  if $T_1,\ldots,T_N$ is irreducible then it is completely reducible with $k=1$. 

We can now formulate our main theorem for the path-dependent shrinking target set defined as in \eqref{def-STS-real}. 
\begin{theorem}\label{main}
Let $\{f_1, \dots, f_N\}$ with $f_i=T_i+a_i$ be an affine iterated function system. Suppose that $\|T_i\|<\tfrac 12$ for all $1\leq i \leq N$,  that $T_1,\ldots,T_N$ is completely reducible, and that every $T_i$ is invertible.  Let $\ell$ be as in Definition \ref{def:path}. Let $\mu$ be any ergodic measure on $(\Sigma, \sigma)$.

Then for $\mu$-almost every choice of $\bj$, the limit defining the pressure \eqref{eq:pressure} exists and is independent of the choice of $\bj$. Further, for these $\bj$, for Lebesgue almost all $a_1, \dots, a_N\in \R^d$, the Hausdorff dimension of the path-dependent shrinking target set $R^*(\bj, \ell)$ is given by $\min\{s_0, d\}$, where $s_0$ is the unique value for which $P(s_0, \bj)=0$. 

\end{theorem}

\begin{remark}
The focus in the current work is on the thermodynamics of shrinking target pressures for weakly quasimultiplicative potentials. From a more geometric perspective, it would also be interesting to study, for some H\"older continuous function ${\psi}: \Lambda \rightarrow \mathbb{R}$, the path-dependent shrinking target set
\[
\left\{x\in \Lambda: E^n x \in B\Big(z_0, \exp\big\{-\sum_{i=0}^{n-1} \psi(E^ix)\big\}\Big), \ \text{for infinitely many } n \right\}.
\]
General resolution of problems of this kind seems currently well out of reach, as the symbolic description of geometric balls is a very involved problem in general: A ball centred at $z_0$ might intersect many cylinders from $\Sigma$, and there is no easy way to determine which ones are intersected and which ones are not. In this sense, geometric balls are incompatible with the dynamical system. 

For some special classes of self-affine sets, the non-path-dependent version of this geometric shrinking target set has been studied \cite{BR, JK}. The structure of the particular classes of self-affine sets handled in these articles allows a translation between geometric and symbolic languages, which is not available in the general case that we handle here. We also point out that in both \cite{BR, JK}, the imposed assumptions force the relevant potentials to be multiplicative, assumption so strong that the thermodynamical issues that we deal with in the current work are removed  entirely. 

\end{remark}

\begin{remark} If we consider an IFS with constant Lyapunov exponents, then our result recovers the Hausdorff dimension of the set $\mathcal{S}(z_0, \psi)$ from \cite{U02}. However, when this is not the case, we investigate here a situation that has not been seen in any past work, including \cite{HV97, BW14, LWWX14, R11}. The difference is, that when we project the targets $\bj(\bi,n)$ to balls in $\Lambda$, the radii of the balls not only depend on $\bi$ ($\pi(\bi)$), but also on $\bj$, or more precisely, the Lyapunov exponent at $\bj$. In this sense Theorem \ref{main} is new even in the conformal case. \end{remark}

The article is organized as follows. In Section \ref{sec:upper} we give a simple proof for the upper bound of the Hausdorff dimension and make other preliminary observations. In Section \ref{sec:modular} we describe the general framework of studying the dynamics in multiplicative blocks, and in particular we investigate various forms of the pressures for weakly multiplicative potentials. This method is likely to be applicable more widely than just in the context of shrinking targets. In Section \ref{sec:targetpressure} we  specialize to singular value potentials relevant in the shrinking target problem. In Section \ref{sec:lowerbound} we apply these to define a Cantor set and a mass distribution that are used to prove the lower bound of the Hausdorff dimension of the shrinking target set in the case of the $\mu$-typical $\bj$.


\section{Preliminaries and the upper bound}\label{sec:upper}

In this section, we first prove some preliminary technical lemmas. Then, we go on to prove that the zero point of the limsup pressure defined above always gives an upper bound to the Hausdorff dimension of the shrinking target set $R^*(\bj, \ell)$.

We have the following two lemmas on the properties of the length function $\ell$. Recall its definition from Definition \ref{def:path}. 
\begin{lemma}\label{lem:l-prop} 
If for all $n\geq 1$ there exists $\bj\in \Sigma$ such that $\ell_n(\bj)<\kappa$, then there exists $\bi$ such that for infinitely many $n$'s we have $\ell_n(\bi)<3 \kappa$.
\end{lemma}
\begin{proof} We first assert that for all $\bi$ and for all $n\in \mathbb{N}$, we have $\ell_n(\bi) \geq -\kappa$. Otherwise, if there exists some $\bi\in \Sigma$, and some $n\in\mathbb{N}$, such that $\underline{\ell}(\bi|_n)=\ell_n(\bi)<-\kappa$, then by \eqref{estimate-ell}, for all $k$,
$$\ell_{kn}((\bi|_n)^\infty)=\underline{\ell}((\bi|_n)^k) \leq \underline{\ell}(\bi|_n)+(k-1)(\underline{\ell}(\bi|_n)+\kappa) \to -\infty \quad (k\to\infty),$$
 which contradicts with our assumption \eqref{assmp:ln-infinity}.
 
 Then, we assert that if for some $\bi\in \Sigma$, and some $n\in\mathbb{N}$, $\ell_n(\bi)<\kappa$, then for all $m<n$ we have $\ell_m(\bi)<3\kappa$. Otherwise, if for some $m<n$,  $\ell_m(\bi)\geq 3\kappa$. Then, by \eqref{estimate-ell}, and the first assertion, we have 
\[
\ell_n(\bi) \geq \ell_m(\bi) + \ell_{n-m}(\sigma^m\bi) -\kappa \geq 3\kappa-\kappa -\kappa=\kappa,
\]
which is a contradiction.
  
 Now, by assumption, for all $n\in \mathbb{N}$, there exists $\bi_n$ such that $\ell_n(\bi_n)<\kappa$.  Since we have only $N$ choices for the first symbol for the infinite sequence $(\bi_n)_{n\geq 1}$, we can find a symbol $a_1$ which appears as the first symbol in infinitely many infinite words $\bi_n$. Hence, by the second assertion, we have $\ul(a_1)<3\kappa$. Similarly, let $a_2$ be the symbol such that infinitely many infinite words $\bi_n$ begin with $a_1a_2$. Then $\ul(a_1 a_2) < 3\kappa$. Continuing with this process, we will obtain an infinite word $\bi=a_1a_2\dots$ such that $\ell_n(\bi)<3\kappa$ for all $n\geq 1$.

\end{proof}

\begin{lemma}\label{lem:prop-l}
There exist some $L_{\max}> L_{\min}>0$ and $\kappa'>0$, such that for all $n \geq 1$
\begin{equation} \label{eqn:bounds}
n L_{\min} - \kappa' \leq \ell_n(\bi) \leq n L_{\max} + \kappa'.
\end{equation}
\end{lemma}
\begin{proof}
By assumption \eqref{assmp:ln-infinity} and Lemma  \ref{lem:l-prop}, there exists $n_0\in \mathbb{N}$, such that $\ell_{n_0}(\bi) > \kappa$ for all $\bi\in\Sigma$.
 Denote by $\kappa_{\min},  \ \kappa_{\max}$ the minimum and maximum of $\{\ell_{n_0}(\bi): \bi\in\Sigma\}$. Denote by $\kappa_{\min}', \ \kappa_{\max}'$ the minimum and maximum of $\{\ell_k(\bi): k\in [0,n_0-1], \ \bi\in\Sigma\}$. Then, one can easily check that for every $\bi\in\Sigma$ and for every $n$ we have
\[
\left\lfloor \frac n {n_0} \right\rfloor (\kappa_{\min}-\kappa ) + \kappa_{\min}' \leq \ell_n(\bi) \leq \left\lfloor \frac n {n_0} \right\rfloor ( \kappa_{\max}+\kappa ) + \kappa_{\max}',
\]
which gives \eqref{eqn:bounds} with $L_{\min}=(\kappa_{\min}-\kappa)/n_0$, $L_{\max}=(\kappa_{\max}+\kappa)/n_0$ and  $\kappa'=\max\{\kappa_{\max}, \ \kappa_{\min}-\kappa'\}$.

\end{proof}

We will prove that the limsup pressure $s\mapsto P^*(s, \bj)$ has a unique zero. The proof relies on the following standard lemma, that we state for completeness.

\begin{lemma}\label{lem:singular}
Let $(T_1,\ldots,T_N)$ be invertible linear maps in $\R^d$. Denote
\[
\alpha_-:=\max_{1\leq j \leq N} \big\{ \|T_j^{-1}\| \big\}, \quad \alpha_+:=\max_{1\leq j \leq N} \big\{ \|T_j\| \big\}. 
\]
 Let $\ri\in \Sigma_n$ be a finite word. Let $0\leq t<s\leq d$. Then
\[
(\alpha_-)^{-n(s-t)} \leq \frac {\phi^s(T_\ri)} {\phi^t(T_\ri)} \leq (\alpha_+)^{n(s-t)}.
\]  
In particular,
\[
(\alpha_-)^{-ns} \leq \phi^s(T_\ri) \leq (\alpha_+)^{ns}.
\]
\end{lemma}
\begin{proof}
Both are straightforward consequences of the definition, by the facts that $\alpha_1(TU)\le \alpha_1(T)\alpha_1(U)$ and $\alpha_d(TU)\ge \alpha_d(T)\alpha_d(U)$ for any two linear maps $T$ and $U$. 
\end{proof}

\begin{lemma}\label{lem:zero}
Fix $\bj\in \Sigma$. There is a unique $s_0$ such that the limsup pressure $P^*(s_0, \bj)=0$.  
\end{lemma}

\begin{proof}
By Lemma \ref{lem:prop-l}, there exist $L_{\max}>L_{\min}>0$ and $\kappa'>0$ satisfying 
\[
L_{\min}|\ri| -\kappa'\le \underline{\ell}(\ri)\le L_{\max}|\ri|+\kappa', \quad \forall \ri\in \Sigma_* \ \text{with} \ |\ri| \to \infty.
\]
Thus, by Lemma \ref{lem:singular}, for any $\bj\in \Sigma$ and for any $k\in \mathbb{N}$, we have
\[
\alpha_-^{-(\kappa'+L_{\max}k)\delta} \le \frac{\sum_{|\ri|=k}\phi^{s+\delta}(\ri\bj|_{\underline{\ell}(\ri)})}{\sum_{|\ri|=k}\phi^s(\ri\bj|_{\underline{\ell}(\ri})}\le \alpha_+^{(\kappa'+L_{\min}k)\delta}, \quad \forall s, \delta>0.
\]
 Then, it follows that $P^*(s, \bj)$ is continuous and strictly decreasing in $s$. Further, $P^*(0, \bj)>0$ and $P^*(s, \bj)\to -\infty$ as $s\to \infty$. Therefore, a unique zero always exists. 
\end{proof}

The following consequence of complete reducibility will be essential to our analysis:
\begin{theorem}\cite[Theorem 5.1]{MorKae}\label{thm:mksingular}
Let $s>0$. If $\{T_1, \dots, T_N\}$ is completely reducible collection of invertible maps
, there is a constant $c>0$, a number $p\le {d\choose {\floor{s}}}{d\choose \ceil{s}}$, and weakly quasimultiplicative potentials $\Psi^s_1, \dots, \Psi^s_p$ taking values in $(0,1)$ such that 
\[
c^{-1}\phi^s(\ri)\le \max_{j=1, \dots, p}\Psi^s_j(\ri)\le c\phi^s(\ri)
\]
and
\[\min_{j=1,\dots,p} \Psi_j^s(\ri) \geq \|T_\ri^{-1}\|^{-s}\]
for all $\ri\in \Sigma$. 
\end{theorem}

The potentials $\Psi_j^s$ have the following further properties: For all $\bi\in \Sigma$, the value $\Psi^{s}_j(\bi|_n)$ is decreasing in $n$. Further, given $s>0$, the potential $\Psi^{s}_j$ evaluates almost exactly the same on two words that differ by just few letters on one end. That is, for $M =\max_i \{\|T_i\|^s, \|T^{-1}_i\|^s\}>0$, for  any pair $\ri, \rj\in \Sigma^*$ such that $\ri|_n=\rj|_n$, we have 
\begin{equation}\label{eq:nottoosmallPot}
\Psi^{s}_j(\ri) \le M^{|\rj|-n} \Psi^{s}_j(\rj). 
\end{equation}
We will need to refer to this property also for general potentials $\Psi$, so we record it here: We say that $\Psi$ is non-collapsing, if for some $A>0$, and all  $\ri, \rj\in \Sigma^*$ such that $\ri|_n=\rj|_n$, it satisfies
\begin{equation}\label{eq:nottoosmall}
\Psi(\ri) \le A^{|\rj|-n} \Psi(\rj). 
\end{equation}

The few following lemmas will need to be applied for more than one potential, so we write them in a general form. To that end, given a potential $\psi$ on $\Sigma^*$, set the notation 
\begin{equation}\label{eq:Zj}
Z(\bj)=Z(\psi, \bj)=\lim_n \tfrac 1n \log \psi(\bj|_n). 
\end{equation}
We will show now that this limit often exists for almost every $\bj$. 

\begin{lemma}\label{lem:Z}
Let $\bj\in \Sigma$ and let $\mu$ be an ergodic measure on $(\Sigma, \sigma)$. Let $\psi$ be any sub-multiplicative potential. Then, for $\mu$-almost every $\bj\in \Sigma$ the limit $Z(\psi, \bj)=\lim\tfrac 1n \log \psi(\bj|_n)$ exists and takes a common value independent of $\bj$. 

\end{lemma}

\begin{proof}
Define 
\[
X(m, n)=\log \psi(\sigma^m(\bj)|_{n-m}).
\]
Then applying Kingman's sub-additive ergodic theorem to $X(m,n)$, we obtain the existence and almost everywhere uniqueness of 
\[
Z(\bj)=\lim_n \tfrac 1n X(0,n). 
\]
\end{proof}

We will need to define many different notions of pressure over the course of the proofs. For now, given a sub-multiplicative potential $\psi$, denote by $\Pfull(\Sigma, \psi)$ the limit
\[
\Pfull(\Sigma, \psi)=\lim_{n\to \infty} \frac 1n \log \sum_{|\ri|=n}\psi(\ri). 
\]

\begin{lemma}\label{lem:pressurexist}
Let $\mu$ be any ergodic measure on $\Sigma$, and let $\psi$ be any sub-multiplicative potential. Then the limit in the definition of $\Pfull(\Sigma, \psi)$ exists. 

\end{lemma}

\begin{proof}
Once again, the statement follows immediately from Kingman sub-additive ergodic theorem applied to 
\[
X(m,n)=\log \sum_{|\ri|=n}\psi (\sigma^m(\ri)). 
\] 
\end{proof}

We will need another notion of pressure for an arbitrary weakly  quasi-multiplicative potential $\Psi$, that is well-aligned with the shrinking target set. Set 
\begin{equation}\label{eq:QMtargetpressure}
\Pfull(\Psi, \bj)=\lim_{n\to \infty} \frac 1n \log \sum_{|\ri|=n} \Psi(\ri \bj|_{\ul(\ri)})
\end{equation}
This limit often exists, as the following lemma demonstrates. 

\begin{lemma}\label{lem:targetpressurexist}
Let $\mu$ be any ergodic measure on $\Sigma$, and let $\Psi$ be any weakly quasi-multiplicative potential, with $\Psi(\rj)$ decreasing in $\rj$ and satisfying \eqref{eq:nottoosmall}. Let $\ul$ satisfy the assumptions of Definition \ref{def:path}, and let $\bj\in \Sigma$, be generic in the sense of Lemma \ref{lem:Z}. Then the shrinking target pressure limit \eqref{eq:QMtargetpressure} exists. 
\end{lemma}
\begin{proof}
Define 
\[
Q_n=\sum_{|\ri|=n}\Psi(\ri \bj|_{\ul(\ri)}). 
\]
Let $\bj$ be generic in the sense of Lemma \ref{lem:Z}, that is, the following limit exists
\[
Z(\Psi, \bj)=\lim_{n\to \infty} \tfrac 1n\log \Psi(\bj|_n). 
\]
Set
\[
\widetilde Q_n=\sum_{|\ri|=n}\Psi(\ri )e^{Z(\Psi, \bj) \cdot\ul(\ri)}. 
\]
It is straightforward to see that $\ri \mapsto \Psi(\ri ) e^{Z(\Psi, \bj)\cdot\ul(\ri)}$ is a sub-multiplicative potential, and hence, by Lemma \ref{lem:pressurexist}, the limit $\lim_{n} \tfrac 1n \log \widetilde Q_n$ exists. 
We will prove that 
\begin{equation}\label{eq:same}
0=\lim_{n\to \infty} \tfrac 1n \log \frac{\widetilde Q_n}{Q_n}, 
\end{equation}
 from which the claim then immediately follows. 
 
In this proof we will use the notation $\sim$ to compare two quantities, when they can be estimated from above and beneath by each other, up to a constant. Firstly, note that for some constant $c$, $Q_n\le c\widetilde Q_n$ by submultiplicativity of $\Psi$, and since $e^{Z(\Psi, \bj)\ul(\ri)}\sim \Psi(\bj|_{\ul(\ri)})$ for large $|\ri|$. Hence, it suffices to check that $Q_n\ge c' \widetilde Q_n$ for some constant $c'$.

To that end, let $n, k\in \N$, and let $|\ra|=k, |\ri|=n$. We are aiming to compare $\Psi(\ri\ra\bj|_{\ul(\ri\ra)})$ to $\Psi(\ri\ra\bj|_{\ul(\ri)})$. Notice that if $\ul(\ri\ra)\ge\ul(\ri)$, by \eqref{eq:nottoosmall}
\[
\Psi(\ri\ra\bj|_{\ul(\ri\ra)}) \ge M^{\ul(\ri\ra)-\ul(\ri)}\Psi(\ri\ra\bj|_{\ul(\ri)}), 
\]
and otherwise 
\[
\Psi(\ri\ra\bj|_{\ul(\ri\ra)}) \ge \Psi(\ri\ra\bj|_{\ul(\ri)}). 
\]
Here, for every pair $|\ra|=k, |\ri|=n$,
\[
|\ul(\ri\ra)-\ul(\ri)|\le |\ul(\ri\ra)-\ul(\ri)-\ul(\ra)|+|\ul(\ra)|\le k L_{\max}+\kappa'+\kappa, 
\]
by Lemma \ref{lem:prop-l} and Definition \ref{def:path}, so that  
\[
\Psi(\ri\ra\bj|_{\ul(\ri\ra)}) \ge M^{k L_{\max}+\kappa'+\kappa}\Psi(\ri\ra\bj|_{\ul(\ri)}). 
\]

Now, denote $\ul(K)=\max_{|\ri|\le K} \ul(\ri)$. Then for some constant $c$ that depends only on $K$, 
\begin{align*}
Q_n&\ge c \sum_{m=n}^{n+K+\ul(K)} Q_m\ge c' \sum_{|\ri|=n}\Psi(\ri\ \con(\ri, \bj|_{\ul(\ri)}) \bj|_{\ul(\ri)}),
\end{align*}
where the second inequality follows from the above, and from the fact that the connecting word for $\ri$ and $\bj_{\ul(\ri)}$ must appear amongst the $Q_m$ for $m=n, \dots, n+K+\ul(K)$. Continuing from here, 
\begin{align*}
&\ge c'\sum_{|\ri|=n}\Psi(\ri)e^{Z(\Psi, \bj)\cdot\ul(\ri)}\\
&=c'\widetilde Q_n. 
\end{align*}

This finishes the proof.

\end{proof}

Now, we are ready to give the upper bound of the Hausdorff dimension.
\begin{lemma}
For every $\bj \in \Sigma$ and every $\ell$ as in Definition \ref{def:path}, the Hausdorff dimension of $R^*(\pi(\bj), \ell)$ is bounded from above by $\min\{s_0, d\}$, where $s_0$ is the unique real number satisfying $P^*(s_0, \bj)=0$.
\end{lemma}
\begin{proof}
This is a standard affine covering argument. We provide the details of the proof for the convenience of the reader. 

The Hausdorff dimension is always bounded from above by $d$. For the upper bound $s_0$, let $s>s_0$ be arbitrary. Notice that $R^*(\pi(\bj))$ is a limsup set, so that for all $n$ it is covered by  
\[
\bigcup _{|\ri|=n}\pi[\ri\bj|_{\underline{\ell}(\ri)}]. 
\]
By the definition of the singular value function, each of the cylinders $\pi[\ri\bj|_{\underline{\ell}(\ri)}]$ can be covered by
$c\phi^s(T_{\ri\bj|_{\underline{\ell}(\ri)}})\alpha_{\lfloor s \rfloor+1}(\ri\bj|_{\underline{\ell}(\ri)})^{-s}$ cubes of sidelength $\alpha_{\lfloor s \rfloor+1}(\ri\bj_{\underline{\ell}(\ri)})$,
where $c$ is an absolute constant.  
Hence, 
\begin{equation}\label{eq:convergence}
\mathcal H^s(R^*(\bj, \ell))\le \lim_{n\to \infty}c\sum_{|\ri|=n} \phi^s(T_{\ri\bj|_{\underline{\ell}(\ri)}}). 
\end{equation}
Since $s>s_0$, we have $\sum_{n=1}^\infty\sum_{|\ri|=n} \phi^s(T_{\ri\bj|_{\underline{\ell}(\ri)}})<\infty$. Therefore, \eqref{eq:convergence} implies that $\mathcal H^s(R^*(\bj, \ell))=0$, completing the proof of the lemma. 
\end{proof}

\section{$R$-modular symbolic spaces and pressure}\label{sec:modular}

To estimate the Hausdorff dimension of the shrinking target set in Section \ref{sec:targetpressure}, we need to consider the pressure formulas for the singular value potential $\log \phi^s$. However, in the interest of more general applicability, in this section we will define and investigate a general weakly quasi-additive potential $\psi$. Hence, let $\psi$ now be a fixed, weakly quasi-additive potential. 

The idea of this section is as follows: We know that $\psi$ is weakly quasi-additive. We can write all sequences of  $\Sigma$ in blocks of some fixed length $R$, and simply insert the necessary connecting words between these blocks -- we call this modified space an $R$-modular space. This change makes the potential essentially multiplicative on this modified space. Of course, we need to check that this procedure does not change any quantities that we care about, in particular, the pressure. The comparison of the pressure defined on the full symbolic space and the $R$-modular space, Corollary \ref{cor:compfull} is the main result of this section.

Recall that $\psi$ being weakly quasi-additive means that there are constants $Q, K$ such that for all $\ri, \rj\in \Sigma_*$ we can find $\rk\in \Sigma_*$, $|\rk|\le K$ such that 
\begin{equation}\label{eq:Q}
\psi(\ri\rk \rj)\ge Q+\psi(\ri)+\psi(\rj). 
\end{equation}
Let 
\begin{equation}\label{eq:L}
L=\sup_{\ri\in \Sigma_*} |\psi(\ri)|/|\ri|.
\end{equation}
\begin{definition}[$R$-modular words]
Fix a large positive integer $R$. We say that a word $\ri\in \Sigma_*$ is $R$-{\it modular} if it can be presented in the form
\[
\ri = \rr_1 \rk_1 \rr_2 \ldots \rk_{n-1} \rr_n,
\]
where the words $\rr_m$ have length $R$ and $\rk_m=\con(\rr_1\rk_1\dots\rk_{m-1}\rr_m, \rr_{m+1} )$ are the connecting words of length $|\rk_m|\le K$ for weak quasi-additivity: for each $m$,
\[
\psi(\rr_1 \rk_1 \ldots \rr_{m+1}) \geq Q + \psi(\rr_1 \rk_1 \ldots \rr_m) + \psi(\rr_{m+1}).
\]
The set of $R$-modular words is denoted by $\mathcal M_R$. The set of $R$-modular words of length $n$ is denoted by $\mathcal M_R^n$. 
\end{definition}

\begin{definition}[$R$-modular extension]
Let $\ri_1\in \Sigma_*$. We say that $\ri_2\in \Sigma_*$ is an $R$-{\it modular extension} of $\ri_1$ if for some $n\in \N$
\[
\ri_2 = \ri_1 \rk_0 \rr_1 \rk_1 \rr_2 \ldots \rk_{n-1} \rr_n,
\]
where for every $m$ we have $|\rr_m|=R$, and $\rk_m=\con(\ri_1\rk_0\rr_1\rk_1\dots\rr_m, \rr_{m+1})$ are the connecting words for weak quasi-additivity. 
\end{definition}

The definition of $R$-modular extension is recursive: modular extension of a modular extension is a modular extension. Also, the $R$-modular words can be thought of as $R$-modular extensions of the empty word (with $\rk_0$ also chosen as the empty word).

We want to apply thermodynamic formalism on $R$-modular spaces, thus we need to define an $R$-{\it modular pressure}. 

\begin{definition}[$R$-modular pressure]
Given a potential $\psi$ defined on $\Sigma_*$, define the $R$-{\it modular pressure} as the limit
\begin{equation}\label{eq:Rmodpre}
\Pmod(\psi) := \limsup_{n\to\infty} \frac 1n \log \sum_{k=0}^{R+K-1} \sum_{\ri\in \mathcal M_R^{n+k}} e^{\psi(\ri)}.
\end{equation}
\end{definition}

\begin{remark}
For a word of length $n$, the lengths of its $R$-modular extensions do not start before $n+R$, and sometimes even $n+R+K$. This is the reason for the sum over $k=0, \dots, R+K-1$ in the definition of $R$-modular pressure. 
\end{remark}

We would like to compare the $R$-modular pressure to the ordinary pressure of the potential $\psi$ on $\Sigma_*$. In what follows, we shall have to vary the space on which we consider the pressure, and hence we use for the most part the following slightly cumbersome notation. 

\begin{definition}[Full pressure]
Given a potential $\psi$ defined on $\Sigma_*$, define the full pressure corresponding to $\psi$ to be 
\begin{equation}\label{eq:fullpre}
\Pfull(\Sigma, \psi)=\lim_{n\to \infty} \frac 1n \log \sum_{|\ri|=n}e^{\psi(\ri)}. 
\end{equation}
\end{definition}

It is immediate from the definitions that 
\[
\Pmod(\psi) \leq \Pfull(\Sigma, \psi).
\]
One of our main goals in this subsection will be showing that the opposite inequality, or at least something close to it, also holds.

Up until now we worked with the symbolic space $\Sigma$ defined on an alphabet $A$ of size $N$, but from here on we will also need to consider the alphabet $A^R$ and the symbolic space $\Sigma^R$ built over it. Given a potential $\psi$ defined on $\Sigma_*$, there is a natural corresponding potential $\tilde{\psi}$ defined on $A^R$ by the formula
\[
\tilde{\psi}(\ra) = \psi(\ra), 
\]
and on longer words in the alphabet $A^R$, extend $\tilde \psi$ by additivity: for $\ra_1,\ldots, \ra_n\in A^R$, let 
\[
\tilde{\psi}(\ra_1 \ra_2 \ldots \ra_n) = \tilde{\psi}(\ra_1) +\ldots +\tilde{\psi}(\ra_n).
\]
Hence $\tilde{\psi}$ coincides with $\psi$ on the words of length $R$, but not necessarily on the words of length $2R$, $3R$, and so on.

For a word $\ra_1\ldots \ra_n \in (A^R)^n$, we have an associated $R$-modular word with blocks $\rr_1=\ra_1,\ldots, \rr_n=\ra_n$. Conversely, for an $R$-modular word $\ri$ with blocks $\rr_1,\ldots, \rr_n$, we can associate  the word $\rr_1 \rr_2 \ldots \rr_n \in (A^R)^n$. Unfortunately, neither of these associations gives us a well defined map. Fortunately, the multivalued maps obtained in this way are not too-multi-valued, as the following lemma shows. 

\begin{lemma}\label{lem:numberof}
Given a word $\ra_1\ldots \ra_n \in (A^R)^n$, there are at most $K_0^{n-1} $ $R$-modular words with blocks $\rr_1=\ra_1,\ldots, \rr_n=\ra_n$, where $K_0:=1+N+\ldots +N^K$. These words have length between $Rn$ and $Rn+K(n-1)$.

Given an $R$-modular word $\ri$ of length $m$, it can be divided into $R$-blocks with legal connecting words in at most $(K+1)^{m/R}$ ways. The number of $R$-blocks in these representations is between $m/(R+K)$ and $m/R$.
\end{lemma}

\begin{proof}

For the first claim: to know an $R$-modular word with prescribed blocks $\rr_1,\ldots,\rr_n$, we need still to choose the connecting parts $\rc_1,\ldots, \rc_{n-1}$, and each of them can be chosen in no more than $K_0:=1+N+\ldots +N^K$ ways. Thus, for a given word $\ra_1\ldots \ra_n \in (A^R)^n$ we get no more than $K_0^{n-1}$ $R$-modular words built with these blocks. Further, the obtained $R$-modular words have length between $Rn$ and $Rn+K(n-1)$.

In the opposite direction, given an $R$-modular word of length $m$, there might be many ways in which it can be divided into blocks. Basically, we need to mark the beginning of each $R$-block, and only then we have the full information. The distance between the consecutive beginnings of the $R$-blocks varies between $R$ and $R+K$, that is no more than $m/R$ times that we need to make a choice and we have at most $K+1$ possibilities each time. Thus, for a given word of length $m$ it can be presented as an $R$-modular word in at most $(K+1)^{m/R}$ ways, and the number of $R$-blocks in these representations is between $m/(R+K)$ and $m/R$.
\end{proof}

We need one more notion of pressure, defined as follows: 
\begin{equation}\label{eq:Rfullpre}
\Pfull(\Sigma^R, \tilde{\psi})=\lim _{n\to \infty} \frac 1n \log \sum_{\ra_1,\ldots, \ra_n\in A^R}e^{\tilde \psi(\ra_1\dots\ra_n)}. 
\end{equation}
As an intermediate step in the comparison of $\Pmod(\psi) $ from \eqref{eq:Rmodpre} and $\Pfull(\Sigma, \psi)$ from \eqref{eq:fullpre}, we will compare $\Pfull(\Sigma^R, \tilde\psi)$ with the pressure $\Pmod(\Sigma, \psi)$. More precisely, we will compare the $n$-level approximation to $\Pfull(\Sigma^R, \tilde{\psi})$ with the $Rn$-level approximation to $\Pmod(\Sigma, \psi)$.

Denote the partial sums
\begin{equation}\label{eq:fullsum}
\Pfulln{n} = \sum_{\ra_1,\ldots, \ra_n\in A^R} e^{\tilde{\psi}(\ra_1 \ldots \ra_n)}
\end{equation}
and
\begin{equation}\label{eq:eq:modsum}
\Pmodn{n} = \sum_{\ri\in \mathcal M^n_R} e^{\psi(\ri)}.
\end{equation}

In the statement of the following definition, and what follows after, it is good to keep in mind that to account for all $R$-modular words, we need to sum over $\mathcal M_R^{n+k}$ for $k=0, \dots, R+K-1$ in the definition of $R$-modular pressure. 

\begin{proposition} \label{prop:compmod}
Let $n, R\in \N$. There exist $c_1, c_2$ depending on $\psi, Q, K$ but not on $R$ such that
\[
\left| \frac 1{Rn} \log \Pfulln{n} - \frac 1 {Rn} \log \sum_{m=0}^{R+K-1} \Pmodn{Rn+m}\right| \leq \frac{c_1}n + \frac{c_2}R.
\]
\end{proposition}

\begin{proof}
Given an $R$-modular word $\ri$ built on $R$-blocks $\rr_1,\ldots, \rr_n$, we have
\begin{equation}\label{eq:closepsi}
Q(n-1) \leq \psi(\ri) - \sum_{k=1}^n \tilde{\psi}(\rr_k) \leq LK(n-1).
\end{equation}
We use this chain of inequalities together with Lemma \ref{lem:numberof} to compare $\Pfulln{n}$ and $\Pmodn{n}$ to each other. There are two immediate consequences of \eqref{eq:closepsi}, which we will state next.  

The first consequence of \eqref{eq:closepsi} follows like this: by Lemma \ref{lem:numberof}, for every summand $e^{\tilde{\psi}(\ra_1\dots \ra_n)}$ appearing in $\Pfulln{n}$, we have at most $K_0^{n-1}$ corresponding summands $e^{\psi(\ri)}$ belonging to some $\Pmodn{m}$  with $m\in [Rn, Rn+K(n-1)]$. Each of these summands is bounded from above by $ e^{\tilde{\psi}(\rr_1)+\ldots +\tilde{\psi}(\rr_n)} \cdot e^{KL(n-1)}$. Thus,
\begin{equation}\label{eq:Pmodbigger}
\sum_{m=Rn}^{Rn+K(n-1)} \Pmodn{m} \geq K_0^{n-1} e^{KL(n-1)} \Pfulln{n}=: C_1^{n-1}\Pfulln{n}.
\end{equation}

The second consequence is deduced as follows: by Lemma \ref{lem:numberof}, for every summand $e^{\psi(\ri)}$ appearing in $\Pmodn{n}$ we have at most $(K+1)^{n/R}$ corresponding summands $e^{\tilde{\psi}(\ra_1\dots \ra_n)}$ belonging to some $\Pfulln{m}$  with $m\in [n/(R+K),1+n/R]$.  By \eqref{eq:closepsi}, each of these summands is at most $Q^{-n/R} e^{\psi(\ri)}$. Thus, 
\begin{equation}\label{eqn:szac1}
\sum_{m=n/{(R+K)}}^{1+R/n} \Pfulln{m} \geq (K+1)^{n/R} Q^{-n/R} \Pmodn{n}=: C_2^{n/R}\Pmodn{n}. 
\end{equation}

That is, if $\Pfulln{n}$ is large then there must be some $m_0\in[Rn, Rn+K(n-1)]$ such that $\Pmodn{m_0}$ is large, and the same in the other direction.

To deduce the statement of the proposition from \eqref{eq:Pmodbigger} and \eqref{eqn:szac1}, we will need an intermediate fact. This is the following inequality, which we will next prove for any $m_0\ge Rn$: 
\begin{equation} \label{eqn:regul2}
 \Pmodn{m_0} \leq {C_3}^{m_0-Rn}\sum_{k=0}^{R+K-1} \Pmodn{Rn+k}, \quad \text{with} \ C_3= Ne^L.
\end{equation}
Indeed, for each $\ri\in A^{m_0}$, and for each $k=0, \dots, m_0-Rn$, there is $\ri' =\ri|_{Rn+k}\in A^{Rn+k}$ which is modular such that
\begin{equation}\label{eq:wordtomod}
e^{\psi(\ri)} \leq e^{\psi(\ri')+(m_0-Rn-k)L}\leq e^{(m_0-Rn)L}e^{\psi(\ri')},
\end{equation}
and there are at most $N^{m_0-Rn-k}$ such words $\ri$ corresponding to one $\ri'$. This implies
\[
\Pmodn{m_0}=\sum_{\ri \in \mathcal{M}_R^{m_0}} e^{\psi(\ri)} \leq \sum_{k=0}^{R+K-1} N^{m_0-Rn-k} \cdot e^{(m_0-Rn)L} \sum_{\ri \in \mathcal{M}_R^{Rn+k}} e^{\psi(\ri)}.
\]
Hence,
\[
\Pmodn{m_0}\leq N^{m_0-Rn} \cdot e^{(m_0-Rn)L} \sum_{k=0}^{R+K-1} \Pmodn{Rn+k},
\]
which proves \eqref{eqn:regul2}. 

We will now combine these estimations to prove the claim. First, note that \eqref{eqn:regul2} implies
\[
\sum_{m=Rn}^{Rn+K(n-1)} \Pmodn{m} \leq \sum_{m=Rn}^{Rn+R+K-1} \Pmodn{m} \cdot \left(1 + \sum_{m=Rn+R+K}^{Rn+K(n-1)} C_3^{m-Rn} \right).
\]
Substituting to \eqref{eq:Pmodbigger}, we get
\[
C_1^{n-1} \Pfulln{n} \leq \left(1 + \sum_{r=R+K}^{K(n-1)} C_3^{r} \right) \sum_{m=Rn}^{Rn+R+K-1} \Pmodn{m}.
\]
Taking logarithms, we have 
\begin{align*}
\log \Pfulln{n}- \log \sum_{m=Rn}^{Rn+R+K-1} \Pmodn{m} \leq  {\log \left(1 + \sum_{r=R+K}^{K(n-1)} C_3^{r} \right)} - {(n-1)\log C_1}.
\end{align*}
This implies one of the inequalities needed to conclude the proposition.

In the opposite direction, we want to carry on from \eqref{eqn:szac1}. To this end, note that $\tilde{\psi}$ is an additive potential and constant on the first level cylinders. Thus, for any $m\in \mathbb{N}$, we have   
\[
\Pfulln{m} = \Pfulln{1}^m,
\]
which implies
\begin{align*}
&\sum_{m=\lceil n/{(R+K) \rceil}}^{\lfloor n/R \rfloor} \Pfulln{m}=\sum_{m=\lceil n/{(R+K) \rceil}}^{\lfloor n/R \rfloor} \Pfulln{1}^m
\end{align*}
Hence,
\begin{align*}
&\sum_{m=\lceil n/{(R+K) \rceil}}^{\lfloor n/R \rfloor} \Pfulln{m}\\
=&\sum_{m=\lceil n/{(R+K) \rceil}}^{\lfloor n/R \rfloor} \Pfulln{1}^{m-\lfloor n/R \rfloor} \cdot \Pfulln{\lfloor n/R \rfloor}.  \\
=& \sum_{k=0}^{\lfloor n/R \rfloor-\lceil n/(R+K) \rceil} \max\left\{1, \ \frac 1 {\Pfulln{1}}\right\}^{k} \cdot \Pfulln{\lfloor n/R \rfloor}\\
\leq& \left(\frac {nK} {R(R+K)}+1 \right)  \max\left\{1, \ \frac 1 {\Pfulln{1}}\right\}^{nK/R(R+K)} \cdot \Pfulln{\lfloor n/R \rfloor}.
\end{align*}
Observe that
\[
\Pfulln{1} \geq N^R e^{-R\sup_{i\in A}|\psi(i)|}=:C_4^R,
\]
Then, combining this with \eqref{eqn:szac1}, we obtain
\[
 \left(\frac {nK} {R(R+K)}+1 \right)   C_4^{nK/(R+K)} \cdot \Pfulln{\lfloor n/R \rfloor} \geq C_2^{n/R} \Pmodn{n}.
\]
Therefore,
\begin{align*}
&\sum_{m=0}^{R+K-1}  \left( \frac {(Rn+m)K} {R(R+K)}+1 \right)   C_4^{(Rn+m)K/(R+K)} \cdot \Pfulln{\lfloor (Rn+m)/R \rfloor} \\
\geq & \sum_{m=0}^{R+K-1} C_2^{(Rn+m)/R} \Pmodn{Rn+m}.
\end{align*}
This implies
\begin{align*}
&(R+K)  \left(\frac {(Rn+R+K-1)K} {R(R+K)}  +1 \right)C_4^{(Rn+R+K-1)K/(R+K)} \cdot \Pfulln{\lfloor(Rn+R+K-1)/R\rfloor} \\
\geq & C_2^n \sum_{m=0}^{R+K-1} \Pmodn{Rn+m}.
\end{align*}
Taking logarithms, we have 
\begin{align*}
&\log \Pfulln{\lfloor(Rn+R+K-1)/R\rfloor}- \log\sum_{m=0}^{R+K-1} \Pmodn{Rn+m} \\
\geq & n \log C_2- \log \left( n+1+{K(K-1) \over R}+R+K\right) -\left({(Rn-1)K \over R+K} +K\right)  \log C_4.
\end{align*}

Finally, remark that when $R\geq K$, 
\begin{align*}
&\Pfulln{\lfloor(Rn+R+K-1)/R\rfloor}\\ \leq&  \Pfulln{n+1}\\
=&\Pfulln{1} \cdot \Pfulln{n}\\
\leq & N^Re^{R\sup_{i\in A}|\psi(i)|}\cdot \Pfulln{n}. 
\end{align*}
This finishes the proof of the second inequality needed for the claim, and hence the proof of the proposition.

%
%

\end{proof}

Let us now look at
\begin{equation}\label{eq:sumpsi}
\Pordn{n}:=\sum_{\ri\in A^n} e^{\psi(\ri)}. 
\end{equation}
We want to compare it with $\Pmodn{n}$, which comparison we will achieve via $\Pfulln{n}$ and Proposition \ref{prop:compmod}.

\begin{corollary} \label{cor:compfull}
There exist $c_3, c_4$ depending on $\psi, Q, K$ but not on $R$ such that
\[
\left| \frac 1n \log \Pordn{n} - \frac 1n \log \sum_{m=0}^{R+K-1} \Pmodn{n+m}\right| \leq \frac {c_3}n + \frac {c_4}R.
\]
\end{corollary}

\begin{proof}
Let $\ri\in\Sigma_*$ be a word of length $Rn$. We can divide $\ri$ into $n$ words of length $R$
\[
\ri=\rr_1 \rr_2 \ldots \rr_n
\]
and we see that
\[
\psi(\ri) \leq \sum_{k=1}^n \tilde{\psi}(\rr_k).
\]
Thus, 
\begin{equation} \label{eqn:full}
\Pordn{Rn} \leq \Pfulln{n}.
\end{equation}

We have already mentioned the trivial fact that $\Pordn{n} \geq \Pmodn{n}$. Thus, the claim is a corollary of Proposition \ref{prop:compmod}. 
\end{proof}

To summarize the results of this subsection: By Corollary \ref{cor:compfull} the $R$-modular sequences, on which the potential is quasi-additive, carry (for large $R$) almost the full pressure. The $R$-modular space of modular sequences is not invariant under the shift map, but it is closely related to $\Sigma^R$. This relation is going to be sufficient to obtain regularity properties needed to estimate the dimension-like characteristics of sets constructed with use of the modular words and modular extensions. We will look at this in the context of shrinking targets in the next section.

\section{$R$-modular shrinking target pressure}\label{sec:targetpressure}

In this section we will apply the $R$-modular pressure analysis of the previous section to the problem of shrinking targets. Consider a completely reducible, invertible affine iterated function system $\{f_1, \dots, f_N\}$, and a length sequence $\ell$ satisfying Definition \ref{def:path} fixed. 

Let $s>0$. Recall that in this completely reducible situation there are weakly quasimultiplicative potentials $\Psi^s_j$ which well approximate $\phi^s$, by Theorem \ref{thm:mksingular}. In other words, for some $c>0$, for all $\ri$,
\begin{equation}
\label{eq:comparemax}
c^{-1}\phi^s(\ri)\le \max_{j=1, \dots, p}\Psi^s_j(\ri)\le c\phi^s(\ri), 
\end{equation}
and there $Q, K$ such that for all $\ri,\rj\in \Sigma_*$ there is some connecting word $\rc=\con(\ri, \rj)\in \Sigma_*$ with $|\rc|\le K$ and 
\begin{equation}\label{eq:weakquas}
Q^{-1} \Psi^s_j(\ri)\Psi^s_j(\rj)\le \Psi^s_j(\ri\rc\rj)\le Q \Psi^s_j(\ri)\Psi^s_j(\rj). 
\end{equation}
It follows from \eqref{eq:comparemax} immediately that 
\[
\lim _{n\to \infty } \frac 1n \log \sum_{|\ri|=n}\phi^s(\ri)=\lim _{n\to \infty } \frac 1n \log \sum_{|\ri|=n}\max_{j=1,\dots, p}\Psi^s_j(\ri). 
\]
Furthermore, by \eqref{eq:weakquas}, we can apply Lemma \ref{lem:pressurexist} to all of $\Psi_j^s$. Let $\bj$ be generic in the sense of Lemma \ref{lem:Z} with respect to an ergodic measure $\mu$ on $\Sigma$, and all $\Psi^s_j$ simultaneously. Recall that the shrinking target pressure corresponding to the potential $\ri\mapsto\log \Psi^s_j(\ri\bj|_{\ul(\ri)}))$ is denoted by $\Pfull(\Psi^s_j, \bj)$, 
\[
\Pfull(\Psi^s_j, \bj)=\lim_{n\to\infty} \frac 1n \log \sum_{|\ri|=n}\Psi^s_j(\ri\bj|_{\ul(\ri)}). 
\]
This limit exists by Lemma \ref{lem:targetpressurexist}. 

We have the following comparison between this pressure and the shrinking target pressure $P(s, \bj)$: 
\begin{lemma}
Let $s>0$, $\bj$ be generic in the sense of Lemma \ref{lem:Z} for all potentials $\Psi_j^s$ from Theorem \ref{thm:mksingular}.  Then 
\[
P(s, \bj)=\max_{j=1,\dots, p} \Pfull(\Psi^s_j, \bj). 
\]
\end{lemma}
\begin{proof}
From \eqref{eq:comparemax} we know that for each $n$
\[
c^{-1}\max_{j=1, \dots, p} \sum_{|\ri|=n} \Psi^s_j(\ri\bj|_{\ul(\ri)}) \le \sum_{|\ri|=n}\phi^s(\ri\bj|_{\ul(\ri)})\le c\sum_{|\ri|=n} \sum_{j=1}^p  \Psi^s_j(\ri\bj|_{\ul(\ri)})\le cp\max_{j=1, \dots, p} \sum_{|\ri|=n}\Psi^s_j(\ri\bj|_{\ul(\ri)}). 
\]
As the potentials $\Psi_j^s$ are weakly quasi-multiplicative, we know by Lemma \ref{lem:targetpressurexist} and the choice of $\bj$ that each of the limits $\Pfull(\Psi_j^s, \bj)$ exists. The claim follows from this immediately. 
\end{proof}

Now, let $s_0$ be the root of the pressure $P(s, \bj)$, i.e. $P(s_0, \bj)=0$. Then, from the above lemma, it is immediate that for some $j=1, \dots, p$, 
\[
\Pfull(\Psi_j^{s_0}, \bj)=0. 
\]
We will now carry on with the analysis for this quasimultiplicative potential.

\subsection{Shrinking target pressure}
Throughout this section, we fix one of the potentials $\Psi_j^{s_0}$ such that $\Pfull(\Psi_j^{s_0}, \bj)=0$. Fix $\bj$ to be such that $Z(\Psi^{s_0}_j, \bj)=\lim_n \tfrac 1n\log\Psi^{s_0}_j(\bj|_n)$ exists. By Lemma \ref{lem:Z} this is a generic property.
Denote $\psi_{s_0, \bj}:\ri\mapsto \log\Psi^{s_0}_j(\ri\bj|_{\ul(\ri)})$. Then, in the notation of the previous section, 
\begin{equation}\label{eq:targetpressure}
\Pfull(\Sigma, \psi_{s_0, \bj})=\Pfull(\Psi^{s_0}_j, \bj). 
\end{equation}
In this subsection we have all the sums over the whole space $\Sigma$ (or the whole $R$-modular space), so to simplify notation we will leave out $\Sigma$ from the notation for the most part. That is, we denote the pressure by
\begin{equation}
\Pfull(\psi_{s_0, \bj})=\Pfull(\Sigma, \psi_{s_0, \bj})=\Pfull(\Psi^{s_0}_j, \bj). 
\end{equation}

The first thing we need to do is to modify the pressure in such a way that it works well on modular words. Given $\ri\in \Sigma_*$, we denote $\l(\ri):= \con(\ri,\bj|_{\ul(\ri)})\bj|_{\ul(\ri\rc)}$. Recall that $\con(\ri, \bj|_{\ul(\ri)})$ is a connecting word of length at most $K$ such that
\[
\Psi_j^{s_0}(\ri \con(\ri, \bj|_{\ul(\ri)})\bj|_{\ul(\ri)}) \geq e^Q \Psi_j^{s_0}(\ri) \Psi_j^{s_0}(\bj|_{\ul(\ri)}).
\]
We note that by Definition \ref{def:path} and Lemmas \ref{lem:l-prop} and \ref{lem:prop-l}, $\bj|_{\ul(\ri\con(\ri, \bj|_{\ul(\ri)}))}$ and $\bj|_{\ul(\ri)}$ differ by at most $2\kappa+KL_{\max}$ letters at the end. Furthermore, the potential $\Psi^{s_0}_j$ evaluates almost exactly the same on two words that differ by just few letters on one end by \eqref{eq:nottoosmall}. Hence, the same $\con(\ri, \bj|_{\ul(\ri)})$ satisfies also 
\[
\Psi^{s_0}_j(\ri\con(\ri, \bj|_{\ul(\ri)})\bj|_{\ul(\ri\con(\ri, \bj|_{\ul(\ri)}))}) \geq Q' \Psi^{s_0}_j(\ri) \Psi^{s_0}_j(\bj|_{\ul(\ri)}).
\]
for some uniform constant $Q'$. We can then define $\psi_{s_0, \l}:\ri\mapsto \log\Psi^{s_0}_j(\ri\l(\ri))$, and it follows that 
\begin{equation}\label{eq:pressurepsipartial}
\Pfull^*(\psi_{s_0, \l}):= \limsup_{n\to\infty} \frac 1n \log \sum_{|\ri|=n} \Psi^{s_0}_j(\ri\l(\ri)) = \limsup_{n\to\infty} \frac 1n \log \sum_{|\ri|=n} \Psi^{s_0}_j(\ri) \Psi^{s_0}_j(\bj|_{\ul(\ri)}).
\end{equation}
Denote by 
\begin{equation}\label{eq:sumtargetpsi}
\Sfull(\psi_{s_0, \bj}, n)=\sum_{|\ri|=n}e^{\psi_{s_0, \bj}(\ri)} \quad \textrm{ and }\quad \Sfull(\psi_{s_0, \l}, n)=\sum_{|\ri|=n}e^{\psi_{s_0, \l}(\ri)}
\end{equation}
the sums at level $n$ in the definitions of $\Pfull^*( \psi_{s_0, \bj})$ and $\Pfull^*(\psi_{s_0, \l})$.

The following lemma shows that the difference between $\psi_{s_0, \bj}$ and $\psi_{s_0, \l}$ is asymptotically unimportant. 
\begin{lemma}\label{lem:comparepressure}
\[
P_{\rm full}^*(\psi_{s_0, \bj})=P_{\rm full}^*(\psi_{s_0, \l}).
\]
\end{lemma}
\begin{proof}
By weak quasi-multiplicativity, for every $\ri\in \Sigma_*, |\ri|=n$ we have, 
\[
\Psi^{s_0}_j(\ri\l(\ri)) \geq Q'\Psi^{s_0}_j(\ri\bj|_{\ul(\ri)}).
\]
Hence,
\[
Q'\Sfull(\psi_{s_0, \bj}, n)\le \Sfull(\psi_{s_0, \l}, n), 
\]
which implies that 
\[
P_{\rm full}^*(\psi_{s_0, \bj})\le P_{\rm full}^*(\psi_{s_0, \l})
\]
On the other hand, the word $\ri\con(\ri,\bj|_{\ul(\ri)})$ has length at most $n+K$, thus, the summand $\Psi^{s_0}_j(\ri\l(\ri))$ appears as one of the terms in $\Sfull(\psi_{s_0, \bj}, n+k), k=0,\ldots,K$. Therefore,
\[
\Sfull(\psi_{s_0, \l}, n) \leq \sum_{k=0}^{K} \Sfull(\psi_{s_0, \bj}, n+k).
\]
The claim follows from here. 
\end{proof} 

Now that we have defined a 'multiplicative' version of the pressure $\Pfull^*(\psi_{s_0, \bj})$, we can create a modular version of it. Fix some positive integer $R$ and remember the definition of $R$-modular words $\mathcal M^R$ from the previous section. Define
\begin{equation}\label{eq:pressuremodpsipartial}
\Pmod^*(s,\bj) := \Pmod^*(\psi_{s_0, \l})=\limsup_{n\to\infty} \frac 1n \log \sum_{k=0}^{K+R-1}\sum_{\ri \in \mathcal M_R^{n+k}} \Psi^{s_0}_j(\ri\l(\ri)).
\end{equation}

 Denote 
\begin{equation}\label{eq:summodpsipartial}
S_{\rm mod}( \psi_{s, \l},  n) = \sum_{\ri\in \mathcal M_R^n} \phi^s(\ri\l(\ri)).
\end{equation}
The main result of this subsection is the following result. 

\begin{proposition}\label{prop:shrink}
There exist $c_5$ and $c_6$ depending on $Q$ and $K$, but not on $R$, such that
\[
\left|\frac 1n \log S_{\rm full}( \psi_{s_0, \l}, n) - \frac 1n \log \sum_{m=0}^{R+K-1} S_{\rm mod}( \psi_{s_0, \l},  n+m)\right| \leq \frac {c_5}n + \frac {c_6}R.
\]
\end{proposition}
\begin{proof}
In Corollary \ref{cor:compfull} we have proved an analogous result for a general $\psi$. To apply Corollary \ref{cor:compfull}  in the shrinking target context, the proof will only require a minimal modification. Let $\ri$ be a word of length $n$ which gives us a summand $\Psi^{s_0}_j(\ri\l(\ri))$ appearing in $\Sfull(\psi_{s_0, \l},  n)$. We divide $\ri$ into blocks of size $R$ and construct an $R$-modular word $\ri'$ out of them by inserting the connecting words. Then, for some uniform constant $c>0$, we have all of the following:
\begin{itemize}
\item[i)] $|\ri'| \leq |\ri| \cdot (1+c/R)$,
\item[ii)] $\Psi^{s_0}_j(\ri') \geq \Psi^{s_0}_j(\ri) \cdot e^{-c|\ri|/R}$,
\item[iii)] $\ul(\ri') \leq \ul(\ri) + c|\ri|/R$, and hence $\Psi^{s_0}_j(\l(\ri')) \geq \Psi^{s_0}_j(\l(\ri)) \cdot e^{-c|\ri|/R}$,
\item[iv)] the map $\ri \to \ri'$ is at most $e^{c|\ri|/R}$-to-1.
\end{itemize}

Indeed, i) and iv) are consequences of Lemma \ref{lem:numberof}, and ii) is the exact counterpart of \eqref{eq:wordtomod} for the potential $\log \Psi^{s_0}_j$. The property iii) follows by combining the approximate additivity property \eqref{estimate-ell} with Lemma \ref{lem:prop-l}, since $\ri'$ is $\ri$ with $|\ri|/R$ inserted additional subwords each of which has length at most $K$.

The points i)-iv) imply the assertion.
\end{proof}

%

\subsection{Shrinking target pressure on abstract modular space $\Sigma^R$} \label{sec:targetmodular}

Fix a large positive integer $R$. Let $\ra\in \Sigma_*$ be a finite word. In this subsection our goal is to describe the set of $R$-modular extensions of the word $\ra$ in the language of alphabet $A^R$. In this subsection, the alphabet varies, so we use a more careful notation again, writing
\[
P(s, \bj)=\Pfull(\Sigma, \psi_{s, \bj})
\]
and so forth. 

Let $s_0$ again be the solution of the equation
\[
\Pfull(\Sigma, \psi_{s_0, \bj})=0, 
\]
and recall that $\psi_{s_0, \bj}(\ri)=\log \Psi^{s_0}_j(\ri)$ for $\ri\in \Sigma_*$. 
As above, we assume that $\log \Psi^{s_0}_j$ is weakly quasi-additive with constants $Q,K$ (see \eqref{eq:weakquasi}). We assume approximate additivity  of $\ul$ as in Definition \ref{def:path}. We also assume that $\bj$ is such that the limit $Z(\Psi^{s_0}_j, \bj)$ from \eqref{eq:Zj} exists. We know by Lemma \ref{lem:Z} that this is a generic property. 

%

Consider a word $\rr_1 \rr_2 \ldots \rr_m \in (A^R)^m$. We can find an $R$-modular extension of $\ra$ of the form $\ra \rk_0 \rr_1 \rk_1 \rr_2 \ldots \rk_{m-1} \rr_m$. As discussed before, generally more than one such extension exists as we might have some freedom in choosing the $\rk_i$'s. To make sure that everything is well-defined, out of these extensions let us choose one, for example let it be the extension for which $\rk_0$ is the first in lexicographical order of all possible $\rk_0$'s, then $\rk_1$ is first in lexicographical order of all possible $\rk_1$'s under the condition that $\rk_0$ is already chosen, and so on. We denote the resulting $R$-modular extension by $\pi_\ra(\rr_1\ldots \rr_m)$.

This gives us a well-defined map $\pi_\ra:\Sigma^R_* \to \Sigma_*$. Moreover, $\pi_\ra(\rr_1 \rr_2 \ldots \rr_{m+1})$ is an $R$-modular extension of $\pi_\ra(\rr_1 \rr_2\ldots \rr_m)$, so $\pi_\ra$ preserves the cylinder structure of $\Sigma_*$. We also have that
\[
R\leq |\pi_\ra(\rr_1 \rr_2 \ldots \rr_{m+1})| - |\pi_\ra(\rr_1 \rr_2 \ldots \rr_m)| \leq R+K.
\]
Hence for every {$\brr\in \Sigma^R$} and every $n>0$ there exists a smallest $z(\ra,\brr,n)$ such that
\[
n \leq |\pi_\ra(\rr_1 \rr_2 \ldots \rr_{z(\ra,\brr,n)})| \leq n+R+K-1,
\]
where $\rr_1 \rr_2 \ldots \rr_{z(\ra,\brr,n)}$ is an initial segment of $\brr$. 
Naturally, for every $\brr$ we have
\[
\frac n {R+K} \leq z(\ra,\brr,n) \leq 1+\frac n R.
\]

We will denote
\begin{equation}\label{eq:newspace}
{\Sigma^R_{\ra,n} := \{\rr_1 \rr_2 \ldots \rr_{z(\ra,\brr,n)}; \ \rr\in \Sigma^R\}}.
\end{equation}
The sets of corresponding cylinders $[\rr_1 \rr_2 \ldots \rr_{z(\ra,\brr,n)}]$ form a disjoint cover of $\Sigma^R$.

Let us now define an additive potential $\Psi$ on $\Sigma^R$, constant on first level cylinders: for $\rb\in A^R$ we set
\[
\Psi(\rb) := \log \Psi^{s_0}_j(\rb) + Z(\Psi_j^{s_0}, \bj) \ul(\rb).
\]
Denote
\begin{equation}\label{eq:sumPsi}
\Sfull(\Sigma^R_{\ra, n}, \Psi) := \sum_{\rr_1\ldots \rr_m \in \Sigma^R_{\ra,n}} e^{\Psi(\rr_1)+\ldots +\Psi(\rr_m)}.
\end{equation}

\begin{proposition} \label{prop:symb}
One can find $c_7$, $c_8$, $c_9$, independent of $\ra$ and $R$, such that
\[
e^{-c_7n/R -c_8 - c_9 n} \leq S_{\rm full}(\Sigma^R_{\ra, n}, \Psi) \leq e^{c_7n/R + c_8 + c_9 n}
\]
Moreover, we can make $c_9$ arbitrarily small while keeping $c_7$ fixed.
\end{proposition}

\begin{proof}
We want to compare the pressure sum for $\Psi$ to the pressure sum for $\psi_{s_0, \l}$, and then make use of the choice of $s_0$. The sum with which we are comparing $\Sfull(\Sigma^R_{\ra, n}, \Psi)$  is the sum 
\begin{equation}\label{eq:sumpsipartial}
\Sfull(\Sigma, \psi_{s_0, \l}, n):= \sum_{|\ri|=n} \Psi^{s_0}_j(\ri\l(\ri)),
\end{equation}
where $\psi_{s_0, \l}(\ri)$ is the notation for the potential $\log\Psi^{s_0}_j(\ri\l(\ri))$ 
as above.

Every word $\ri\in \Sigma_*$ with $|\ri|=kR=:n$ can be divided into subwords of length $R$ as $\ri=\rr_1 \rr_2 \ldots \rr_k$, where, $\rr_1, \rr_2, \ldots, \rr_k \in A^R$. Then
\begin{equation}\label{eq:uptoconst}
\log \Psi^{s_0}_j(\ri) \le \log \Psi^{s_0}_j(\rr_1) + \ldots + \log \Psi^{s_0}_j(\rr_k) \le \log \Psi_j^{s_0}(\pi_\emptyset(\rr_1\dots \rr_k))+ O(n/R), 
\end{equation}
where the word $\pi_\emptyset(\rr_1\dots \rr_k)$ has a length bounded from below by $kR$ and above by $k(R+K)$. In particular, undoing the logarithms and recalling that $\Psi^{s_0}_j(\ri)$ is decreasing in length of $\ri$,
\[
\left |\sum_{|\ri|=n=kR} \Psi_j^{s_0}(\ri) - \sum_{|\rr_1|,\dots, |\rr_k|=R }\Psi_j^{s_0}(\rr_1)\dots \Psi_j^{s_0}(\rr_k) \right|\le e^{O(n/R)+o(n)}.
\]
Further, 
\[
\ul(\ri) = \ul(\rr_1) + \ldots + \ul(\rr_k) + O(n/R).
\]
Moreover, we know from Lemma \ref{lem:l-prop} that $\ul(\ri) < L_{\max} n+\kappa'$ and hence
\[
|\log\Psi^{s_0}_j(\l(\ri)) - Z(\bj) \ul(\ri)| = o(n).
\]
Comparing $\Sfull(\Sigma^R_{\ra, n}, \Psi)$ to $\Sfull(\Sigma, \psi_{s_0, \l}, n)$ term by term, we obtain from the above estimates that
\begin{equation}\label{eq:Psiclose}
|\Sfull(\Sigma^R_{\ra, n}, \Psi)-\Sfull(\Sigma, \psi_{s_0, \partial}, n)|=e^{(O(n/R)+o(n))}. 
\end{equation}
Recall from Lemma \ref{lem:comparepressure} that $\Pfull^*(\psi_{s, \l})=\Pfull^*(\psi_{s, \bj})$. By our assumptions, $\Pfull^*(\psi_{s_0, \bj})=0$ and so 
\begin{equation}\label{eq:psitargetsmall}
|\Sfull(\Sigma, \psi_{s_0, \l}, n)|=e^{o(n)}. 
\end{equation}
Here, we  need that the limit defining $\Pfull(\Psi^{s_0}_j, \bj)$ exists, which we know from Lemma \ref{lem:targetpressurexist}.

Combining \eqref{eq:Psiclose} and \eqref{eq:psitargetsmall} gives
\[
|\log \Sfull(\Sigma^R_{\ra, n}, \Psi)-{0}| = O(n/R) + o(n),
\]
as desired. 
\end{proof}

We finish this subsection by defining a measure on $\Sigma^R$, which we will utilize in the next subsection when looking for a measure supported on $R(\bj, \ell)$ that will give the lower bound for the Hausdorff dimension of $R^*(\bj, \ell)$. 

\begin{definition}[Measure on $\Sigma^R$]\label{def:nuR}
Recall the potential $\Psi$ satisfying for $\rb\in A^R$ 
\[
\Psi(\rb) = \log \Psi^{s_0}_j(\rb) + Z(\bj) \ul(\rb),
\]
and extended to $(A^R)^n$ by additivity. Then, for each $\rb\in A^R$, set 
\[
\nu_{R}[\rb] := \frac {\exp(\Psi(\rb))} {w_R}, 
\]
where
\[
w_R=\sum_{\rb_1\in A^R} \exp(\Psi(\rb_1))
\]
is a normalizing factor to guarantee that $\nu_R$ is a unit mass distributed on cylinders. The set function $\nu_R$ extends to a measure on $\Sigma^R$ in the natural way by taking for $\rb_1,\dots,\rb_m\in A^R$
\[
\nu_{R}[\rb_1,\dots,\rb_m] = \nu_R[\rb_1]\cdots\nu_R[\rb_m], 
\]
and hence to a measure on $\Sigma^R$ by the Carath\'eodory's extension theorem. 
\end{definition}

\begin{lemma}\label{lem:normaliser}
There exists a uniform constant $c_{10}$ such that we have a bound for the normalizing factor $w_R$, valid for every $R$:
\[
c_{10}^{-1} \leq 
w_R=\sum_{\rb_1\in A^R} \exp(\Psi(\rb_1)) \leq c_{10}.
\]
\end{lemma}

\begin{proof}
As $\nu_R$ is a probability measure and $\Sigma^R_{\ra, n}$ gives a disjoint cover of $\Sigma^R$, for every $n$ we have
\[
\sum_{\rr_1\ldots \rr_k \in \Sigma^R_{\ra,n}} \nu_{R}[\rr_1\ldots \rr_k] =1.
\]
As for each such $\rr_1\ldots \rr_k\in \Sigma^R_{\ra, n}$, we have that $n/R \geq k \geq n/(R+K)$, we get
\[
\log \Sfull(\Sigma^R_{\ra, n}, \Psi)   = w_R \cdot (n/R \cdot (1+O(1/R))), 
\]
and the assertion follows from Proposition \ref{prop:symb}, after we choose $n$ large enough.
\end{proof}

\section{Lower bound for Hausdorff dimension}\label{sec:lowerbound}

We can now begin our study of the path-dependent shrinking target set $R(\bj,\ell)$. The proof strategy is as follows. We will fix some large $R$ and then find a Cantor subset $W_R\subset R(\bj,\ell)$, constructed with the help of $R$-modular extensions. On the subset $W_R$ we will then distribute a measure $\mu$ such that for any cylinder $[\ri]$ we have 
\begin{equation}\label{eq:measure}
\mu[\ri] \leq c\psi(\ri)
\end{equation} 
for some constant $c>0$, where $\psi$ is some suitably chosen potential. The proof of the lower bound is is an application of the following well-known theorem of Falconer and Solomyak. The theorem is stated here in an altered form (for closed subsets $A\subset \Sigma$), but it follows from the same proofs line by line.  

\begin{theorem}[Lemma 3.1 of \cite{F88}, Proposition 3.1 of \cite{S98}]\label{thm:FS}
Consider an affine iterated function system $\{f_1, \dots, f_N\}$ with $f_i=T_i+a_i$, and its corresponding sequence space $\Sigma$ as in Subsection \ref{Subsec:2.3}. Assume that $\|T_i\|<\tfrac 12$ for all $i=1, \dots, N$. 
Let $\mu$ be a finite measure supported on a closed subset $A\subset \Sigma$ such that for all $\rmq\in A, n\in \N$, we have $\mu[\rmq|_n]\le c\phi^s(\rmq|_n)$. Then, for all $t<s$  
\[
\iint |x-y|^{-t}\,d\pi_*\mu\,d\pi_*\mu<\infty,
\]
and in particular $\dim_H\pi(A)\ge s$, for Lebesgue almost all choices of $(a_1, \dots, a_N)$. 
\end{theorem}

\subsection{Construction of the Cantor subset and the corresponding mass distribution} Let us begin with the construction of $W_R$. We fix some fast increasing sequence $(n_i)$, satisfying $n_1 \gg R$ and $n_{i+1}/n_i\to \infty$. Recall from \eqref{eq:newspace} the notation 
\[
\Sigma^R_{\ra, n_1}=\big\{\rr_1\dots\rr_{z(\ra, \rr, n_1)}\mid \rr\in \Sigma^R\big\},
\]
and the definition of $\pi_{\ra}$ from Subsection \ref{sec:targetmodular}. 

For the first step of the construction, we define
\[
W_{R,1} := \bigcup_{\rb_1\in \Sigma^R_{\emptyset, n_1}}  \{\pi_\emptyset(\rb_1) \l(\pi_\emptyset(\rb_1))\}.
\]
Recall the notation 
\[
\Sigma^R_{\emptyset, n_1}=\Big\{\rr_1\dots\rr_{z(\emptyset, \rr, n_1)}\mid \rr\in \Sigma^R\Big\}
\]
and in particular, the union in the definition of $W_{R, 1}$ is taken over finite words of the form $(A^R)^m$ with $m$ varying. The set $W_{R, 1}$ is a collection of finite words, each of the form $\ra \l(\ra)$. 

The second step: 
\[
W_{R,2} := \bigcup_{\ra_1\in W_{R,1}} \bigcup_{\rb_2\in \Sigma^R_{\ra_1, n_2}} \big\{\pi_{\ra_1}(\rb_2) \l(\pi_{\ra_1}(\rb_2))\big\}. 
\]
That is, $W_{R,2}$ consists of words of the form $\ra\l(\ra)$, where $\ra$ are $R$-modular extensions of the words from $W_{R,1}$, such that their length is as close to $n_2$ as possible. And so we carry on:
\[
W_{R,k} := \bigcup_{\ra_{k-1}\in W_{R,k-1}} \bigcup_{\rb_k\in \Sigma^R_{\ra_{k-1}, n_k}} \big\{\pi_{\ra_{k-1}}(\rb_k) \l(\pi_{\ra_{k-1}}(\rb_k))\big\}.
\]
Finally we define $$W_R=\bigcap_{k=1}^\infty W_{R, k}.$$

Every infinite word in $W_R$ is described by a sequence of finite words $\{ \rb_1, \rb_2, \ldots\}$, each $\rb_k$ belonging to the set of $R$-modular extensions of some word determined by the previous $\rb_1,\ldots, \rb_{k-1}$. As we can see, for every point in $W_R$ there are infinitely many times when the initial segment of this word are of the form $\ra\l(\ra)$, hence indeed $W_R \subset R(\bj,\ell)$.

One important observation: for an infinite word in $W_R$ the sequence $(\rb_1, \rb_2,\ldots)$ is in general not uniquely defined. In fact, even in the first step it can happen that two different $\rb_1\in \Sigma^R_{\emptyset, n_1}$ produce the same $\ra_1\in W_{R, 1}$ and the same for the other $\rb_k$ in the sequence. This technicality is important now, as we start distributing a measure on $W_R$.

The construction of the measure $\mu$ is based on the measure $\nu_R$ from Definition \ref{def:nuR}, and goes as follows. First, on cylinders from $W_{R,1}$ we distribute the measure 
\[
\mu_{R,1} := \sum_{\rb_1\in \Sigma^R_{\emptyset, n_1}}  (\pi_\emptyset)_*(\nu_{R}|_{\rb_1}).
\]
That is, for each cylinder $[\pi_\emptyset(\rb_1) \l(\pi_\emptyset(\rb_1))]$ from $W_{R,1}$, we assign the mass $\nu_{R}(\rb_1)$. Note that if there is $\rb_1'\in \Sigma^R_{\emptyset, n_1}$ such that
$$\pi_\emptyset(\rb_1) \l(\pi_\emptyset(\rb_1))=\pi_\emptyset(\rb'_1) \l(\pi_\emptyset(\rb'_1))$$
then the mass of $[\pi_\emptyset(\rb_1) \l(\pi_\emptyset(\rb_1))]$ will be the sum. 


Next, for each $\rb_1\in \Sigma^R_{\emptyset, n_1}$ corresponding to $\ra_1=\pi_\emptyset(\rb_1) \l(\pi_\emptyset(\rb_1))\in W_{R,1}$, we subdivide the measure according to:
\[
\mu_{R,2}|_{[\ra_1]} := \mu_{R, 1}[\ra_1]\sum_{\rb_2\in \Sigma^R_{\ra_1, n_2}} (\pi_{\ra_1})_* (\nu_{R}|_{\rb_2}). \]
That is, for each sub-cylinder $[\pi_{\ra_1}(\rb_2) \l(\pi_{\ra_1}(\rb_2))]$ of $[\ra_1]$, we assign the mass $\mu_{R, 1}(\ra_1) \cdot \nu_{R}(\rb_2)$, with multiplicity if there is repetition. Notice that if it were possible to disregard repetition, we'd have the simple product formula 
\[
\mu_{R, 2}[\pi_{\ra_1}(\rb_2) \l(\pi_{\ra_1}(\rb_2))]=\nu_R(\rb_1)\nu_R(\rb_2). 
\]

In general, assume that some $\ra_1\in W_{R, 1}, \dots, \ra_{k-1}\in W_{R,k-1}$ and the corresponding $\rb_1\in \Sigma^R_{\emptyset, n_1}, \dots, \rb_{k-1}\in \Sigma^R_{\ra_{k-2}, n_{k-1}}$ have been inductively chosen. Then, we define $\mu_{R, k}$ on the cylinders of $W_{R, k}$, by setting for each $\ra_{k-1}\in W_{R, k-1}$, 
\[
\mu_{R,k}|_{[\ra_{k-1}]} = \mu_{R,1}[\ra_1]\cdots\mu_{R, k-1}[\ra_{k-1}]\sum_{\rb_k\in \Sigma^R_{\ra_{k-1}, n_k}} (\pi_{\ra_{k-1}})_*(\nu_R|_{\rb_k}). 
\]
That is, each sub-cylinder $[\pi_{\ra_{k-1}}(\rb_k)\l(\pi_{\ra_{k-1}}(\rb_k))]$ of $[\ra_{k-1}]$ gets assigned the weight 
\[
\mu_{R,1}[\ra_1]\cdots\mu_{R, k-1}[\ra_{k-1}]\nu_R [\rb_k]
\]
with multiplicity if there is repetition. Notice that, again, ignoring repetition would lead to the simple product formula
\begin{equation}\label{eq:product}
\mu_{R,k}[\pi_{\ra_{k-1}}(\rb_k)\l(\pi_{\ra_{k-1}}(\rb_k))]=\nu_{R}[\rb_1]\cdots\nu_{R}[\rb_{k-1}]\nu_R [\rb_k].
\end{equation}
Finally, we take $\mu_R$ as the weak limit of $\mu_{R,k}$.

Let $\rw\in W_R$ and let $\rb_1, \rb_2,\ldots$ be one of its generating symbolic sequences. Let $\ra_1(\rb_1)=\pi_\emptyset(\rb_1) \l(\pi_\emptyset(\rb_1))$, $\ra_2(\rb_1, \rb_2) = \pi_{\ra_1} (\rb_2)  \l(\pi_{\ra_1}(\rb_2)),\ldots$,  $\ra_k(\rb_1, \rb_2, \dots, \rb_k), \dots $ be as above. For any $n$ we want to compare $\mu_R([\rw|_n])$ with $\phi^{s_0}(\rw|_n)$. Let us start by looking at the part of the measure $\mu_R$ coming from the sequence $\rb_1, \rb_2,\ldots$, that is, let us for the time being follow the simplified formula \eqref{eq:product} above where repetition is ignored, and let us denote the consequent product measure by $\tilde \mu$. In that case, for $|\ra_{r-1}|<n\leq |\ra_r|$,
\[
\tilde{\mu}([\rw|_n]) := \nu_{R}([\rb_1]) \cdot \nu_{R}([\rb_2]) \cdot \ldots\cdot \nu_{R}([\rb_{r-1}]) \cdot \nu_{R}([\sigma^{|\ra_{r-1}|}(\rw|_n)]).
\] 

Now that the measure has been defined, we check that it is well-distributed. We begin with a proof that the product version $\tilde\mu$ is. 

\begin{lemma} \label{lem:est}
Let $\rw\in W_R$, $n\in \N$. There exist universal constants $c_{11}, c_{12}, c_{13}$ independent of $n, R, \rw$ such that
\[
\frac {\tilde{\mu}([\rw|_n])} {\phi^{s_0}(\rw|_n)} \leq e^{c_{11} n/R + c_{12} + c_{13} n}.
\]
Moreover, $c_{13}$ can be chosen arbitrarily small keeping $c_{11}$ constant and increasing $c_{12}$.
\end{lemma}
\begin{proof}
We know that $\phi^{s_0}(\rw|_n)=\max_{j=1, \dots, p} \Psi^{s_0}_j(\rw|_{n})$ and hence we can focus on the potential $\Psi^{s_0}_j$. 

Consider first $n\leq n_1$. Let $\rb_1=\rr_1 \rr_2 \ldots \rr_m$, with $\rr_k\in A^R$. The beginning of the sequence $\rw$ is $\pi_\emptyset(\rr_1 \rr_2 \ldots \rr_m) \l(\pi_\emptyset(\rr_1 \rr_2 \ldots \rr_m))$, with $|\pi_\emptyset(\rr_1 \rr_2 \ldots \rr_m)|\approx n_1$.

Recall that $\Psi(\rr_i) = \log \Psi^{s_0}_j(\rr_i) + Z \ul(\rr_i) $.  For $n<n_1$ we have 
\[
\tilde{\mu}([\rw|_n]) = \exp(\sum_{i=1}^k \Psi(\rr_i)+O(k)), 
\]
where $n/(R+K) \leq k\leq n/R$. We also have
\[
\Psi^{s_0}_j(\rw|_n) = \prod_{i=1}^k \Psi^{s_0}_j(\rr_i) \cdot e^{-O(k)}.
\]

By \eqref{eqn:bounds}, we have 
\[
\ul(\ri)\ge \Lmin |\ri|-\kappa', \quad \forall \ri\in \Sigma_*. 
\]
Hence, $\ul(\rr_i) \geq \Lmin R - \kappa'$. Hence, we see that for $n\le n_1$, 
\[
\frac {\tilde{\mu}([\rw|_n])} {\Psi^{s_0}_j(\rw|_n)} \leq e^{O(n/R) +\Lmin Z n}.
\]
Therefore, the assertion follows with a safe margin (we can choose $c_{11}=c_{12}=0$ and also we have an additional exponentially decreasing factor).

This safe margin is immediately used to deal with the case $n_1 < n \leq |\ra_1|$. In this range $\tilde{\mu}([\rw|_n])$ stays constant while $\Psi^{s_0}_j(\rw|_n)$ eventually decreases, by \eqref{eq:comparemax}, as the same is true for $\phi^{s_0}(\rw|_n)$. Hence, we only need to check the situation for $n=|\ra_1|$. There, we have, up to constants,
\[
 \Psi^{s_0}_j(\rw|_{|\ra_1|}) =  \Psi^{s_0}_j(\rw|_{n_1}) \cdot  \Psi^{s_0}_j(\bj|_{\ul(\rw|_{n_1})}).
\]
However, by Lemma \ref{lem:normaliser}, for some $c>0$,
\[
\tilde{\mu}([\rw|_{|\ra_1|}]) = e^{cm} \cdot \prod_{i=1}^m e^{\Psi(\rr_i)}  =  \Psi^{s_0}_j(\rw|_{n_1}) \cdot e^{Z(\bj) \ul(\rw|_{n_1})} \cdot e^{O(n/R)},
\]
hence
\[
\frac {\tilde{\mu}([\rw|_{|\ra_1|}])} { \Psi^{s_0}_j(\rw|_{|\ra_1|})} \leq e^{O(n/R)} \cdot \frac {e^{Z(\Psi^{s_0}_j, \bj) \ul(\rw|_{n_1})}} { \Psi^{s_0}_j(\bj|_{\ul(\rw|_{n_1})})},
\]
and the last factor is sub-exponential.

The proof then follows by induction. If we know the assertion holds for $n=|\ra_{r-1}|$, we write $\rb_r=\rr_1\ldots \rr_m$, and use the formulas
\[
 \Psi^{s_0}_j(\rw|_n) =  \Psi^{s_0}_j(\rw|_{|\ra_{r-1}|}) \cdot \prod_{i=1}^k  \Psi^{s_0}_j(\rr_i) \cdot e^{O(k)}
\]
and
\[
\tilde{\mu}([\rw|_n]) = \tilde{\mu}([\rw|_{|\ra_{r-1}|}]) \cdot \exp(\sum_{i=1}^k \Psi(\rr_i)+O(k)).
\]
Then, almost identical calculations as above give us the assertion for all $|\ra_{r-1}|<n\leq |\ra_r|$.
\end{proof}

Lemma \ref{lem:est} has the following corollary.

\begin{corollary} 
For a measure $\mu_R$ defined using a large enough $R$ in the definition of the $R$-modular space, there is a constant $C>0$ such that for any $\varepsilon>0$,  all $\rw\in W_R, n\in \N$,
\[
\mu_R([\rw|_n]) \leq C\phi^{s_0}(\rw|_n) e^{n\varepsilon}.
\]
\end{corollary}

\begin{proof}
Fix $\varepsilon>0$. We may again focus on the potential $\Psi^{s_0}_j$, as $\phi^{s_0}$ is an upper bound to it. By Lemma \ref{lem:est} we only need to estimate the impact of the overlaps on the measure $\mu_R$. To this end, we need to estimate in how many possible ways a given $R$-modular word (or a given $R$-modular extension) can be obtained. This is a calculation we have already done in Lemma \ref{lem:numberof}: the representation of a word as an $R$-modular word (or $R$-modular extension) is uniquely determined by marking the beginnings of the $R$-blocks, thus a word of length $n$ has at most $(K+1)^{n/R}$ possible representations. Thus, for every $\rw\in W_R$ and every $n$ there are at most $(K+1)^{n/R}$ words $\rw_i\in W_R$ such that
\[
\mu_R[\rw|_n] = \sum_{i=1 }^{(K+1)^{n/R}}\tilde{\mu}[\rw_i|_n].
\]
Together with Lemma \ref{lem:est}, this implies
\[
\mu_R[\rw|_n]\le \exp(c_{11}n/R+c_{12}+c_{13}n) \Psi^{s_0}_j(\rw|_n), 
\]
where we can choose $c_{13}$ arbitrarily small by taking $c_{12}$ larger.
Choosing $R$ large and $c_{13}$ small 
we obtain the claim.
\end{proof}

We are now ready to apply Theorem \ref{thm:FS} and finish the proof of the lower bound as explained at the beginning of the section: Notice that for any $s<s_0$, by Lemma \ref{lem:singular}, 
\[
\mu_R([\rw|_{n}])\le C \phi^s(\rw|_{n})\alpha_+^{n(s_0-s)}e^{n\varepsilon}. 
\]
Choosing $\varepsilon $ small enough, we obtain 
\[
\mu_R([\rw|_{n}])\le C' \phi^s(\rw|_{n}), 
\]
which allows the application of Theorem \ref{thm:FS}. 

\section{Glossary of pressures and partial sums}\label{sec:glossary}

We list in the following tables all the pressure-related notation from the article. 

\begin{table}[h]
\begin{tabular}{l|l|l}
 Notation & In text & Notes \\
 \cline{1-3}
 $\Pmod(\psi) = \limsup\limits_{n\to\infty} \frac 1n \log \sum\limits_{k=0}^{R+K-1} \sum_{\ri\in \mathcal M_R^{n+k}} e^{\psi(\ri)}$ & \eqref{eq:Rmodpre} &   \\
 $\Pfull(\Sigma, \psi)=\lim\limits_{n\to \infty} \frac 1n \log \sum\limits_{|\ri|=n}e^{\psi(\ri)}$ & \eqref{eq:fullpre} &\\
 $\Pfull(\Sigma^R, \tilde{\psi})=\lim\limits_{n\to \infty} \frac 1n \log \sum\limits_{\ra_1,\ldots, \ra_n\in A^R}e^{\tilde \psi(\ra_1\dots\ra_n)}$ & \eqref{eq:Rfullpre} & \makecell{$\tilde \psi$ extended \\from $\psi$ by additivity}\\
\end{tabular}
\caption{Pressures defined from general potential $\psi$.}
\end{table}

\begin{table}[h]
\begin{tabular}{l|l|l}
 Notation & In text & Notes \\
 \cline{1-3}
$\Pmodn{n} = \sum\limits_{\ri\in \mathcal M^n_R} e^{\psi(\ri)}$ & \eqref{eq:eq:modsum} &  \makecell{summed up over $n+k, k=0, \dots, R+K$\\ to obtain pressure} \\
 $\Pordn{n}:=\sum\limits_{\ri\in A^n} e^{\psi(\ri)}$ & \eqref{eq:sumpsi} &\\
 $\Pfulln{n} = \sum\limits_{\ra_1,\ldots, \ra_n\in A^R} e^{\tilde{\psi}(\ra_1 \ldots \ra_n)}$ & \eqref{eq:fullsum}
 & \\
\end{tabular}
\caption{Partial sums for pressures from general potential $\psi$.}
\end{table}

\begin{table}[h]
\begin{tabular}{l|l|l}
 Notation & In text & Notes \\
 \cline{1-3}
$P(s, \bj)=\Pfull(\psi_{s, \bj})=\Pfull(\Sigma, \psi_{s, \bj})$ & \eqref{eq:targetpressure} &  \makecell{$\psi_{s, \bj}:\ri\mapsto \log\phi^s(\ri\bj|_{\ul(\ri)})$ } \\
 $\Pfull^*(\psi_{s, \l})= \limsup\limits_{n\to\infty} \frac 1n \log \sum\limits_{|\ri|=n} \phi^s(\ri\l(\ri))$ & \eqref{eq:pressurepsipartial} & \makecell{$\psi_{s, \l}:\ri\mapsto \log\phi^s(\ri\l(\ri))$,\\ $\l(\ri)=\con(\ri, \bj|_{\ul(\ri)})\bj|_{\underline{\ell}(\ri \con(\ri, \ell|_{\ul(\ri)})})$}\\
$\Pmod^*(\psi_{s, \l})=\limsup\limits_{n\to\infty} \frac 1n \log \sum\limits_{k=0}^{K+R-1}\sum\limits_{\mathcal M_R^{n+k}} \phi^s(\ri\l(\ri))$  & \eqref{eq:pressuremodpsipartial}
 & $\Pmod^*(s,\bj) := \Pmod^*(\psi_{s, \l})$\\
\end{tabular}
\caption{Pressures from the potentials corresponding to shrinking targets.}
\end{table}

\begin{table}[h]
\begin{tabular}{l|l|l}
 Notation & In text & Notes \\
 \cline{1-3}
$\Sfull(\psi_{s, \bj}, n)=\sum\limits_{|\ri|=n}e^{\psi_{s, \bj}(\ri)}$ & \eqref{eq:sumtargetpsi} &  \\
 $\Sfull(\psi_{s, \l}, n)=\sum\limits_{|\ri|=n}e^{\psi_{s, \l}(\ri)}$ & \eqref{eq:sumtargetpsi} &\\
 $S_{\text{mod}}( \psi_{s, \l},  n) = \sum\limits_{\ri\in \mathcal M_R^n} \phi^s(\ri\l(\ri))$ & \eqref{eq:summodpsipartial}
 & \\
 $\Sfull(\Sigma^R_{\ra, n}, \Psi) = \sum\limits_{\rr_1\ldots \rr_m \in \Sigma^R_{\ra,n}} e^{\Psi(\rr_1)+\ldots +\Psi(\rr_m)}$ & \eqref{eq:sumPsi}
 & $\Psi(\rb)=\log \psi^{s_0}(\rb)+Z(\bj)\underline{\ell}(\rb)$\\
\end{tabular}
\caption{Partial sums for shrinking target pressures.}
\end{table}

\newpage

\section*{Acknowledgements}  

The work was partially supported by the France-Poland bilateral project PHC Polonium (44851YC, PPN/BFR/2019/1/00013). M. R. was also partially supported by National Science Centre grant 2019/33/B/ST1/00275 (Poland). HK and MR acknowledge the support of a London Mathematical Society Scheme 2 grant that funded MR's productive visit to Bristol.

\end{document}